\documentclass{article}

\usepackage[utf8]{inputenc}
\usepackage{amsmath}
\usepackage{fullpage}
\usepackage{relsize}
\usepackage{amssymb}
\usepackage{amsfonts}
\usepackage{graphicx}
\usepackage{tikz}
\usepackage{float}
\usepackage{caption}
\usepackage{authblk}
\usepackage[english]{babel}
\usepackage{mathtools}
\usepackage{array}
\usepackage{pifont}
\usepackage{subfig}
\usepackage{xcolor,soul}
\usetikzlibrary[plotmarks]

\makeatletter
\newcommand{\thickhline}{%
    \noalign {\ifnum 0=`}\fi \hrule height 1.5pt
    \futurelet \reserved@a \@xhline
}
\newcolumntype{"}{@{\hskip\tabcolsep\vrule width 1.5pt\hskip\tabcolsep}}
\newcolumntype{?}{!{\vrule width 1.5pt}}
\makeatother

\usepackage[square,comma,numbers,sort&compress]{natbib}
\usepackage{hyperref}

\usetikzlibrary{arrows.meta}
\usetikzlibrary{calc}
\usetikzlibrary{decorations.pathmorphing, patterns,shapes}

\renewcommand{\Re}{\mathrm{Re}}

\renewcommand{\i}{\mathrm{i}}
\newcommand{\e}{\mathrm{e}}
\newcommand{\eps}{\epsilon}
\newcommand{\as}{\mathrm{as}}

\newcommand{\pdiff}[2]{\frac{\partial #1}{\partial #2}}

\renewcommand{\d}{\,\mathrm{d}}
\newcommand{\diff}[2]{\frac{\mathrm{d} #1}{\mathrm{d} #2}}

\title{title}
\author[1]{people\footnote{Electronic address: christopher.lustri@sydney.edu.au}}
\affil[1]{Department of Mathematics and Statistics, 12 Wally's Walk, Macquarie University, New South Wales 2109, Australia}
\date{}

\begin{document}


\begin{center}
\begin{LARGE}
{\sc Locating complex singularities of Burgers' equation using exponential asymptotics and transseries}
\end{LARGE}

\vspace{2ex}
\begin{large}Christopher J.~Lustri$^1$, In\^{e}s Aniceto$^2$, Daniel J.~VandenHeuvel$^3$ \& Scott W.~McCue$^3$\end{large}

\vspace{2ex}
$^1$
School of Mathematics and Statistics, University of Sydney, Sydney NSW 2006, Australia

$^2$
School of Mathematical Sciences, University of Southampton, Southampton SO17 1BJ, United Kingdom

$^3$
School of Mathematical Sciences, Queensland University of Technology, Brisbane QLD 4001, Australia

\end{center}

\vspace{1ex}\noindent Keywords: complex-plane singularities; exponential asymptotics; transseries, resurgence

\vspace{2ex}\noindent {\bf Abstract:} Burgers’ equation is an important mathematical model used to study gas dynamics and traffic flow, among many other applications. Previous analysis of solutions to Burgers' equation shows an infinite stream of simple poles born at $t = 0^+$, emerging rapidly from the singularities of the initial condition, that drive the evolution of the solution for $t > 0$.
We build on this work by applying exponential asymptotics and transseries methodology to an ordinary differential equation that governs the small-time behaviour in order to derive asymptotic descriptions of these poles and associated zeros.
Our analysis reveals that subdominant exponentials appear in the solution across Stokes curves; these exponentials become the same size as the leading order terms in the asymptotic expansion along anti-Stokes curves, which is where the poles and zeros are located.  In this region of the complex plane, we write a transseries approximation consisting of nested series expansions. By reversing the summation order in a process known as transasymptotic summation, we study the solution as the exponentials grow, and approximate the pole and zero location to any required asymptotic accuracy.
We present the asymptotic methods in a systematic fashion that should be applicable to other nonlinear differential equations.

\section{Introduction}

The field of exponential asymptotics, including the associated techniques of analysing so-called beyond-all-orders terms in asymptotic expansions, has a variety of applications in physics and mathematics \cite{BerryHowls1990,chapmanmortimer_2005,Chapman2006,Chapman2009,Chapman_shock_caustic,deng2023exponential,lustri2012,Lustri}, in particular to understand Stokes phenomenon \cite{Berry1988,Chapman,Howls_2012}.  The more advanced use of transseries, transasymptotic summation and the theory of resurgence is less common in applied mathematics \cite{Aniceto_2021,costin_1998,Howls2010bothsides,Howls_2012,oldedaalhuis2005}, although becoming more popular in theoretical physics (see e.g. \cite{Bajnok:2022xgx,Borinsky:2021hnd,Dorigoni:2022bcx,Heller:2021yjh,Jankowski:2023fdz,Reis:2022tni}, as well as \cite{ANICETO20191} and references therein).  Here we apply exponential asymptotics and transseries to locate singularities and zeros of a solution to a second-order ordinary differential equation (ODE) that comes from studying Burgers' equation in a small-time limit, using the method of transasymptotics developed in~\cite{costin1999correlation,costin2002formation,costin2015tronquee}.  Our analysis is motivated in part due to significant contemporary interest in the tracking of complex-plane singularities of solutions to differential equations \cite{Aniceto:2022dnm,Aniceto-upcoming_painleve,Costin_2019,fornberg_2014,OldeDaalhuis_2022,weideman2022}, but also because the methodology presented here should be applicable to other nonlinear ODEs, including those without exact solutions or the kind of special properties that Burgers' equation has (including integrability).

As mentioned, the application we are concerned with is Burgers' equation 
\begin{align}
\pdiff{u}{t} + u \pdiff{u}{x} &= \mu \pdiff{^2 u}{x^2}, \qquad x \in \mathbb{R},\, t > 0,
\label{e:burger}
\end{align}
which is an extremely well-studied nonlinear partial differential equation (PDE). 
On the real line, it illustrates advection-driven wave steepening competing with linear diffusion, and also acts as an idealised model for a number of physical processes including gas dynamics or traffic flow. This PDE was studied in the regime $\mu \ll 1$ using exponential asymptotics by Chapman et al.~\cite{Chapman_shock_caustic}, where the authors described how singularities in the analytic continuation of the PDE solution generate behaviour known as ``Stokes' phenomenon''. By studying Stokes' phenomenon in the solution, the authors explained the onset of smoothed shock fronts at a catastrophe point. Stokes' phenomenon also plays an important role in explaining the small-time behaviour that we consider in the present study.

We consider the $\mu = \mathcal{O}(1)$ regime of Burgers' equation \eqref{e:burger}, supplemented with an initial condition that has simple poles in the complex $x$-plane, appearing as complex-conjugate pairs.  As studied in VandenHeuvel et al.~\cite{vandenheuvel2022burgers} in some detail (as well as Appendix B of \cite{Chapman_shock_caustic}), in the small-time limit there is an inner problem near one of these poles governed by
\begin{equation}\label{eq.z0ODE}
-\frac{1}{2}U - \frac{\xi}{2}\diff{U}{\xi} + U \diff{U}{\xi} = \mu \diff{^2U}{\xi^2},
\end{equation}
where $\xi$ is a similarity variable, subject to the far-field condition
\begin{equation}\label{eq.farfield}
U \sim -\frac{\mathrm{i}}{2\xi} \qquad \mathrm{as} \qquad \xi \to -\i \infty
\end{equation}
(which turns out to act as two boundary conditions). The condition \eqref{eq.farfield} is obtained by matching the solution of \eqref{eq.z0ODE} to a small-time expansion on the real line in $x$. The main objective of our study is to use exponential asymptotics and transseries to estimate the location of the singularities and zeros of $U(\xi)$ in the complex $\xi$-plane.  To interpret the results in terms of the original model, these poles and zeros in the $\xi$-plane are poles and zeros of Burgers' equation 
(\ref{e:burger}) that emerge from the singularities of the initial condition $u(x,0)$ with direction $\arg(\xi)$ and speed $\mathcal{O}(t^{-1/2})$ as $t\rightarrow 0^+$.

The problem (\ref{eq.z0ODE})--(\ref{eq.farfield}) provides an instructive example to apply our methodology for the following reasons.  First, it serves as a prototype for inner problems that arise by treating PDEs in a small-time limit, or similarity solutions more generally.  Second, all singularities of Burgers' equation are simple poles, so our analysis of (\ref{eq.z0ODE})--(\ref{eq.farfield}) avoids any complications from branch cuts and multiple Riemann sheets.  And finally, there is an exact solution of (\ref{eq.z0ODE})--(\ref{eq.farfield}) (which can be written in terms of parabolic cylinder functions; see (\ref{eq:exactsoln})) that can be used to check the results, although this solution itself is difficult to resolve numerically.  While our analysis turns out to be rather detailed in parts, we shall attempt to make the broader strategy accessible to those interested in other differential equations.

We provide a comprehensive background on exponential asymptotics and transseries in Appendix~\ref{sec:background}. Additionally, a brief stripped-back summary of our methodology is as follows:
\begin{itemize}
    \item By direct substitution into  (\ref{eq.z0ODE}), expand the solution $U(\xi)$ as a divergent series
    \begin{equation}
    U\sim \sum_{m=1}^\infty \frac{a_m^{(0)}}{\xi^{2m-1}}
    \quad\mbox{as}\quad\xi\rightarrow\infty,
    \label{eq:Unaiveansatz}
    \end{equation}
    which applies in a sector of the $\xi$-plane including the negative imaginary axis.  This step is straightforward (see section~\ref{sec:formulation}).
    \item Apply techniques in exponential asymptotics to locate the Stokes curves (also known in the literature as Stokes lines) and determine the leading-order form of the exponentially small correction, which turns out to be
    $$
    U_{\mathrm{exp}}\sim 2\pi\mathrm{i}\Lambda\xi^{-\mathrm{i}/2\mu}\mathrm{e}^{-\chi},
    $$
    where $\chi=\xi^2/4\mu$ and $\Lambda$ is a constant (see section~\ref{sec.Stokes}).  The important Stokes curves lie where $\chi$ is real and positive \cite{Dingle}.  This exponentially small term $U_{\mathrm{exp}}$ is said to ``switched on'' as the Stokes curves are crossed in the complex $\xi$-plane.
    \item Extend the series (\ref{eq:Unaiveansatz}) to include full asymptotic expansions for the exponential correction and its exponential correction, and so on, to obtain a transseries
    \begin{equation}
        U \sim \sum_{m=1}^\infty \frac{a_m^{(0)}}{\xi^{2m-1}}+\tau\sum_{m=0}^\infty \frac{a_m^{(1)}}{\xi^{2m-1}}+\tau^2\sum_{m=0}^\infty \frac{a_m^{(2)}}{\xi^{2m-1}}+\ldots 
    = \sum_{n=0}^\infty \tau^n \sum_{m=0}^\infty \frac{a_m^{(n)}}{\xi^{2m-1}}
    \label{eq:Uextendedansatz}
    \end{equation}
    as $\xi\rightarrow\infty$, 
    where $\tau=\sigma\xi^{-(1+\mathrm{i}/2\mu)}\mathrm{e}^{-\chi}$ and $a_0^{(0)}=0$; such an expansion remains well-ordered provided $|\tau|\ll 1$, i.e. away from the anti-Stokes curves (also knowns as anti-Stokes lines) $\Re(\chi)=0$ (see section~\ref{sec:formulatetrans}).
    \item Write the transseries (\ref{eq:Uextendedansatz}) as a transasymptotic summation
    \begin{equation}
    U\sim \sum_{m=0}^\infty \frac{1}{\xi^{2m-1}}\sum_{n=0}^\infty \tau^n a_m^{(n)} =
    \sum_{m=0}^\infty \frac{A_m(\tau)}{\xi^{2m-1}},
    \label{eq:Utransansatz}
    \end{equation}
    where now $\tau$ is no longer required to be small, and therefore (\ref{eq:Utransansatz}) holds in neighbourhood of anti-Stokes curves (see section~\ref{sec:transanalysis})
    \item By direct substitution, determine $A_m$ (section~\ref{sec:transanalysis}) and subsequently locate the singularities (section~\ref{sec:polelocations}) and zeros (section~\ref{sec:zerolocations}) of $A_0(\tau)$, $A_1(\tau)$, etc., which ultimately provides an asymptotic approximation for the singularities and zeros of the original function $U(\xi)$.
\end{itemize}
As mentioned already, the above procedure does not rely on the exact solution of (\ref{eq.z0ODE})--(\ref{eq.farfield}) and can be adopted to other nonlinear ODE problems \cite{costin1999correlation}.  We discuss our results in section~\ref{sec.results} and conclude in section~\ref{sec:conclude}.

\section{Problem Formulation}\label{sec:formulation}

In this section we briefly present the formulation of our problem (\ref{eq.z0ODE})--(\ref{eq.farfield}), summarising relevant parts of Ref.~\cite{vandenheuvel2022burgers}. We first derive the small-time asymptotic behaviour of Burgers' equation (\ref{e:burger}) with an initial condition that is singular at $x = \pm \i$, and observe that the small-time expression ceases to be asymptotic in a small neighbourhood about each of the initial singularities. We then rescale (\ref{e:burger}) to study the behaviour in the neighbourhood of $x=\mathrm{i}$, and find that it is governed by (\ref{eq.z0ODE})--(\ref{eq.farfield})  in a scaled similarity variable.  Finally, for reference, we note the exact solution of (\ref{eq.z0ODE})--(\ref{eq.farfield}) in terms of parabolic cylinder functions.

\subsection{Small-time limit of Burgers' equation}

The background for (\ref{eq.z0ODE})--(\ref{eq.farfield}) starts with Burgers' equation (\ref{e:burger}) with the initial condition
\begin{equation}\label{e:IC}
    u(x,0) = \frac{1}{1+x^2}, \qquad x \in \mathbb{R}.
\end{equation}
This initial condition was chosen as it has a pair of simple poles at $x = \pm \i$, which allows us to more easily apply matched asymptotic expansions to study the small-time behaviour, especially the birth of infinitely many simple poles that are born at $x = \pm \i$ at $t=0^+$ and propagate in the complex $x$-plane for $t>0$~\cite{vandenheuvel2022burgers,Chapman_shock_caustic} (for initial conditions with other types of singularities, including branch points, the analysis is much more complicated~\cite{vandenheuvel2022burgers}).  

By applying a naive outer expansion of the form
\begin{equation}\label{eq.smalltimeseries}
u(x,t) \sim u_0(x) + t u_1(x) + t^2 u_2(x) + \ldots \quad \mathrm{as} \quad t \to 0^+,
\end{equation}
we find by direct substitution that 
\begin{align*}
u_0 &= \frac{1}{1+x^2},\qquad 
u_1 = \frac{2(-\mu + x + 3 \mu x^2)}{(1+x^2)^3},\\
u_2 &= \frac{60\mu^2 x^4 - 32 \mu x + 48 \mu x^3 - x^2(120\mu^2-7) + 12\mu^2 - 1}{(1+x^2)^5}.
\end{align*}
In the limit that $x \to \i$, we have
\begin{equation}
u_0\sim -\frac{\i}{2}\frac{1}{x - \i}, \qquad
u_1\sim  \left(-\frac{1}{4} - \mu \i\right)\frac{1}{(x-\i)^3}\qquad
u_2\sim  \left(\frac{1}{4} - \frac{5\mu}{2} - 6 \i \mu^2\right)\frac{1}{(x-\i)^5}
\label{eq:innerofouter}
\end{equation}
(with similar expressions in the limit $x \to -\i$). Clearly the strongest singular behaviour for $u_0$, $u_1$, $u_2$, $\ldots$, increases in strength by two at each order, which suggests that (\ref{eq.smalltimeseries}) remains valid for $x \in \mathbb{C}$, except for a region in the neighbourhood of the singular points $x = \pm \i$.

\subsection{Inner expansion near $x=\mathrm{i}$}

By comparing terms $u_0$ with $tu_1$ (or $tu_1$ with $t^2u_2$, etc.) near $x=\pm\mathrm{i}$, we see the expansion (\ref{eq.smalltimeseries}) ceases to be well-ordered when $x\mp\i = \mathcal{O}(t^{1/2})$ as $t \to 0^+$. Understanding what happens near the points $x = \pm \i$ is important, because these points generate an infinite number of singularities in the solution $u(x,t)$ of Burgers' equation (\ref{e:burger}), which emerge from the singularities of $u_0(x)$ rapidly for small time.  For what follows, we will concentrate on the singular point $x = \i$, but our analysis applies in the same fashion to $x =-\i$.

The above argument suggests an inner region near $x=\mathrm{i}$, which, in terms of the new coordinates
\begin{equation}\label{eq.innervariables}
\xi = \frac{x - \mathrm{i}}{t^{1/2}},\qquad u = \frac{1}{t^{1/2}}\bar{U}(\xi,t),
\end{equation}
holds for $\xi = \mathcal{O}(1)$ as $t \to 0^+$. In terms of \eqref{eq.innervariables}, Burgers' equation \eqref{e:burger} is written exactly as
\begin{equation}\label{e.intermediate}
    t\pdiff{\bar{U}}{t} - \frac{1}{2}\bar{U} - \frac{1}{2}\xi \pdiff{\bar{U}}{\xi} + U \pdiff{\bar{U}}{\xi} = \mu \pdiff{^2 \bar{U}}{\xi^2}.
\end{equation}
To leading order, we write $\bar{U} \sim U(\xi)$ as $t \to 0^+$, so after substituting into \eqref{e.intermediate} we end up with our main ODE (\ref{eq.z0ODE}).  To match back onto the outer region (\ref{eq.smalltimeseries}), we rewrite (\ref{eq:innerofouter}) in terms of the inner variables (\ref{eq.innervariables}), which gives
\begin{equation}
U\sim -\frac{\mathrm{i}}{2\xi} + \left(-\frac{1}{4} - \i\mu\right)\frac{1}{\xi^3} +  \left(\frac{\i}{4} - \frac{5\mu}{2} - 6 \i \mu^2 \right)\frac{1}{\xi^5} + \mathcal{O}(\xi^{-7})\quad \mathrm{as} \quad \xi \to -\mathrm{i}\infty.
\label{eq:3termfarfield}
\end{equation}
To leading order, this condition reduces to (\ref{eq.farfield}).  As demonstrated in Ref.~\cite{vandenheuvel2022burgers}, a Liouville-Green (WKB) analysis of (\ref{eq.z0ODE})--(\ref{eq.farfield}) shows that only the leading-term of (\ref{eq:3termfarfield}) (i.e., (\ref{eq.farfield})) is required to uniquely determine the solution; the correction terms in (\ref{eq:3termfarfield}) shall be required for our transseries analysis later on (see (\ref{eq:Unaiveansatz}) and (\ref{eq.F0ser})).

\subsection{Exact solution of (\ref{eq.z0ODE})--(\ref{eq.farfield})}\label{sec:exact}

As discussed in Ref.~\cite{vandenheuvel2022burgers}, our main problem (\ref{eq.z0ODE})--(\ref{eq.farfield}) has an exact solution (which, as we have emphasised, we use only for comparing with our asymptotic results), namely
\begin{equation}
    U=\frac{1}{2\sqrt{2\mu}}
    \frac{U_{\mathrm{p}}\left(\frac{1}{2}-\frac{\mathrm{i}}{4\mu},\frac{\mathrm{i}\xi}{(2\mu)^{1/2}}\right)}
    {U_{\mathrm{p}}\left(-\frac{1}{2}-\frac{\mathrm{i}}{4\mu},\frac{\mathrm{i}\xi}{(2\mu)^{1/2}}\right)},
    \label{eq:exactsoln}
    \end{equation}
where $U_{\mathrm{p}}$ is a parabolic cylinder function.  To derive (\ref{eq:exactsoln}), we integrate (\ref{eq.z0ODE}) once to give a Riccati equation, which can be transformed into Kummer's equation via a change of variable.  The solution follows from carefully enforcing the far-field condition (\ref{eq.farfield}).  

To evaluate (\ref{eq:exactsoln}) numerically we can rewrite the solution in terms of hypergeometric functions which we compute using the \texttt{HypergeometricFunctions.jl} package in \textsc{Julia} \cite{hypergeometricfunctions.jl}.  To achieve satisfactory resolution, we have had to compute with arbitrary precision using \textsc{Julia}'s \texttt{ArbNumerics.jl} package \cite{Johansson2017arb,arbnumerics.jl}.  For full details of the derivation and treatment of (\ref{eq:exactsoln}), see Ref.~\cite{vandenheuvel2022burgers}.

\begin{figure}
\centering
\subfloat[$\mu = 0.5$]{
\includegraphics[width=0.45\textwidth]{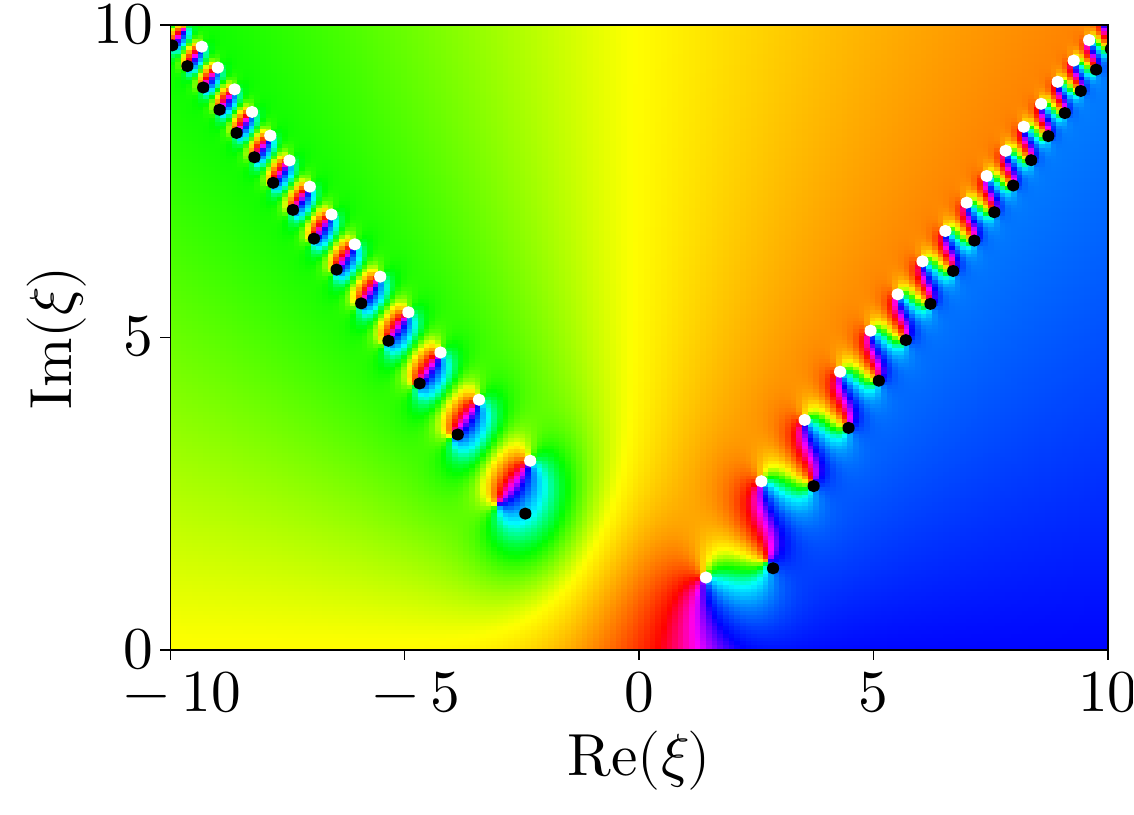}
}
\subfloat[$\mu = 1$]{
\includegraphics[width=0.45\textwidth]{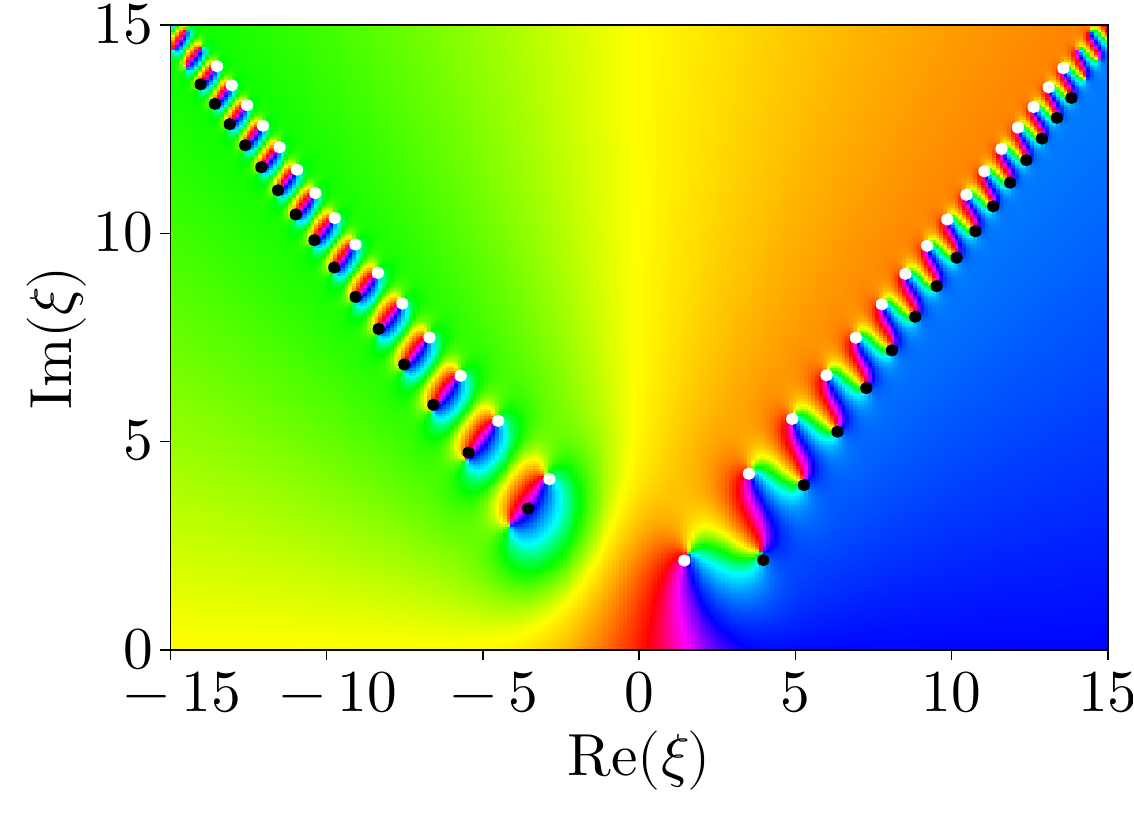}
}

\subfloat[$\mu = 2$]{
\includegraphics[width=0.45\textwidth]{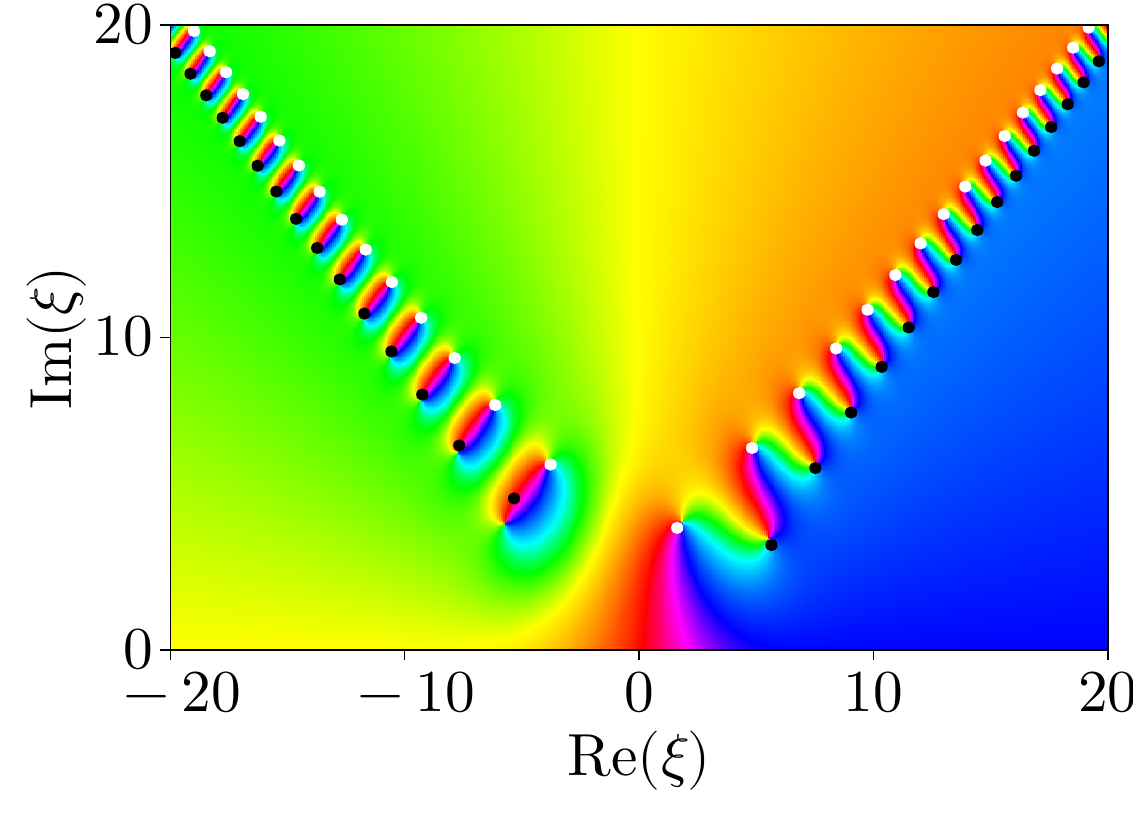}
}
\subfloat[Phase colour indicator]{
\includegraphics[width=0.45\textwidth]{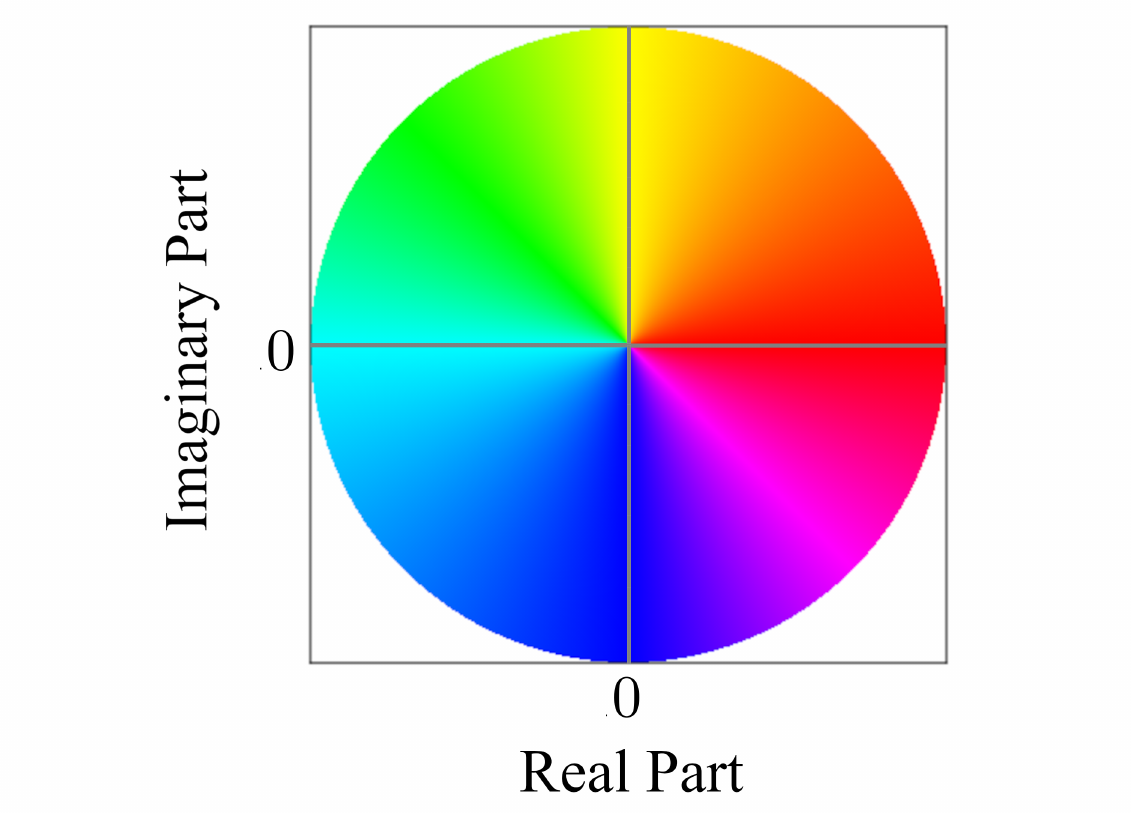}
}
\caption{Phase plots of the solution of (\ref{eq.z0ODE})--(\ref{eq.farfield}) (upper-half plane only) for (a) $\mu=0.5$, (b) $\mu=1$ and (c) $\mu=2$, computed using the exact formula (\ref{eq:exactsoln}).  The colour indicates the phase according to the colour wheel in (d), where for example red represents real and positive.  Also included in (a)-(c) are our asymptotic predictions for pole locations (white dots) via (\ref{e:poleloc}) and zero locations (black dots) via (\ref{e:zeroloc}).}\label{Fig:complex_comp}
\end{figure}

The exact solution (\ref{eq:exactsoln}) is illustrated in Figure \ref{Fig:complex_comp}(a)-(c) via phase portraits for three representative values of $\mu$.  Here we include only the upper-half $\xi$-plane, noting that for these values of $\mu$, there are no poles or zeros in the lower-half plane.  In (d), we include a colour wheel that indicates the phase of the solution~\cite{wegert2010}.  One way to interpret this colour where is to note that locally, a simple zero in (a)--(c) appears locally like the wheel in (d) up to rotation, while a simple pole looks locally like the wheel in (d) with the angular direction reversed.  Importantly, from this figure we can visualise a string of poles and zeros that tend to align themselves at angles $\mathrm{arg}(\xi)=\pi/4$ and $3\pi/4$.  We shall return to this figure later when we discuss the asymptotic predictions of these pole and zero locations.

\section{Stokes switching analysis}\label{sec.Stokes}

\begin{figure} 
\centering
\includegraphics[width=0.85\textwidth]{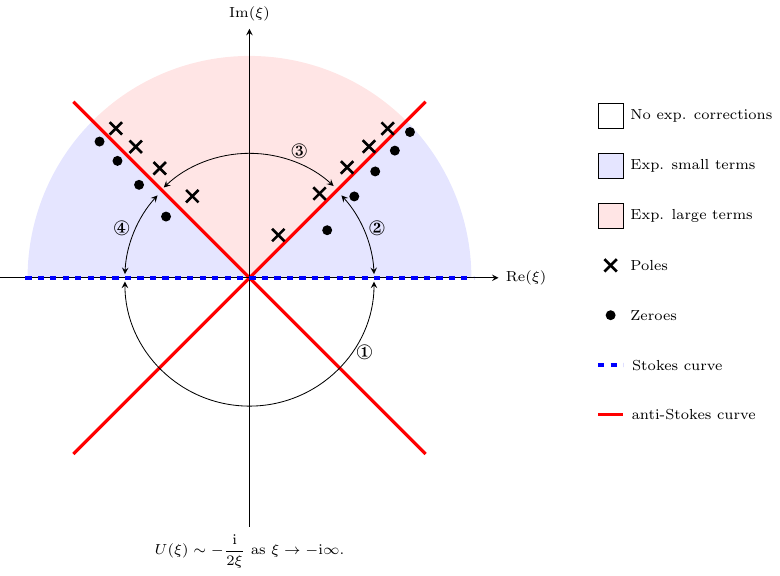}
\caption{Stokes sectors for the inner problem (\ref{eq.z0ODE})--(\ref{eq.farfield}). There are no exponential terms present in the sector containing the negative imaginary axis, where the boundary condition (\ref{eq.farfield}) is specified. Any exponentially small terms that could potentially appear in the expansion would be larger than the boundary behaviour on the negative imaginary axis. This region is shown as Sector \ding{172}. The anti-Stokes curves in the lower half-plane therefore have no effect on the solution. Exponentially small terms appear across the Stokes curves that follow the real axis, in sectors \ding{173} and \ding{175}. As the anti-Stokes curves in the upper half-plane are crossed into Sector \ding{174}, the exponential terms become large, and the asymptotic series terms re-order. The competing exponentials give rise to simple poles and zeros in the solution, represented by crosses and discs, respectively. The poles and zeroes in this figure actual examples computed for $\mu = 1$.}\label{fig.zstokes}
\end{figure}

In the limit $|\xi| \to \infty$, the solution to (\ref{eq.z0ODE})--(\ref{eq.farfield}) exhibits Stokes' phenomenon in the complex $\xi$-plane. This means that subdominant exponential terms are present in different regions of the solution, and that these terms appear as particular curves, known as Stokes curves, are crossed. In this section, we will locate the Stokes curves, which are shown in Figure \ref{fig.zstokes}, and determine the form of the exponential terms. Further, we will determine the location of anti-Stokes curves, or curves across which the exponential terms transition in size from exponentially small to exponentially large.  In the neighbourhood of these anti-Stokes curves, an algebraic power series representation of the solution ceases to be valid, which will lead to a reordering of terms as part of our transasymptotic analysis in section~(\ref{sec.trans}).

\subsection{Power series representation}

To study the asymptotic solution to (\ref{eq.z0ODE})--(\ref{eq.farfield}), we start by setting $\xi = z/\epsilon$ and $U = \epsilon V$, and treat the solution when $z = \mathcal{O}(1)$. In these coordinates, studying the large-$|\xi|$ behaviour is equivalent to taking the small-$\epsilon$ limit while keeping $z = \mathcal{O}(1)$. The exponential asymptotics summarised in this section could be performed in terms of $\xi$, but it is more straightforward in the scaled coordinate $z$; since $\eps$ is an artificial small parameter that is used to index the size of $|\xi|$, it will not appear in the final result \eqref{eq.Uexp} when written in terms of $\xi$. 

In terms of $z$ and $\epsilon$, our governing ODE (\ref{eq.z0ODE}) becomes 
\begin{equation}\label{eq.zeODE}
2\epsilon^2\left(  \mu \diff{^2V}{z^2} - V \diff{V}{z}\right) + z \diff{V}{z} + V = 0.
\end{equation}
Motivated by the far-field condition \eqref{eq.farfield}, we apply the algebraic power series ansatz
\begin{equation}\label{eq.Vser}
V \sim \sum_{n = 0}^{\infty} \epsilon^{2n}V_n(z) \quad \mathrm{as} \quad \epsilon \to 0,
\end{equation}
with 
\begin{equation}\label{eq.V0}
V_0 = -\frac{\mathrm{i}}{2z},
\end{equation}
which describes the solution in a region of the complex plane containing the negative imaginary axis. From the Stokes structure of the solution, we will later be able to determine the broader sector in which this algebraic series is valid. By substituting \eqref{eq.Vser}-(\ref{eq.V0}) into \eqref{eq.zeODE} and matching powers of $\epsilon$, we can obtain a recurrence relation for $n \geq 1$, 
\begin{align}\label{V.recur}
z \diff{V_n}{z} + V_n = - 2 \mu \diff{^2 V_{n-1}}{z^2} + 2\sum_{j=0}^{n-1}\diff{V_j}{z}V_{n-1-j},\qquad
V_n(z) = {o}\left(\frac{1}{|z|}\right) \quad \mathrm{as} \quad z \to -\i\infty.
\end{align} 
The first few terms of the series are therefore given by
\begin{equation}\label{eq.powerseries}
V(z) \sim -\frac{\i}{2z} + \left(\frac{ z-2 \i \mu}{2 z^3}\right)\epsilon^2 + \left( \frac{3 \i z^2+ 16 \mu z-36 \i \mu^2  }{6 z^5}\right)\epsilon^4 + \ldots \quad \mathrm{as} \quad \epsilon \to 0.
\end{equation}
Note that since $V_0$ is singular at $z=0$, then so will $V_1$, $V_2$, $\ldots$.  Now we must determine the region in which this power series approximation is valid, and how the solution behaves outside of this region. We therefore need to understand the Stokes' phenomenon present in $V(z)$.

\subsection{Late-order terms}

A key step in our analysis is determining the asymptotic behaviour of $V_n$ in the limit that $n \to \infty$. The presence of the summation term in (\ref{V.recur}) means that $V_n$ cannot easily be evaluated exactly for arbitrary $n$. Instead, to determine the asymptotic behaviour of the late-order terms, we follow the method proposed by \cite{Chapman} and applied in numerous studies (eg. \cite{chapmanmortimer_2005,Chapman2006,Chapman2009,Chapman_shock_caustic,deng2023exponential,lustri2012,Lustri}) and assume that these terms diverge in a factorial-over-power fashion according to the ansatz
\begin{equation}
    V_n\sim\frac{G(z)\Gamma(n+\gamma)}{\chi(z)^{n+\gamma}}\quad \mathrm{as} \quad n\to\infty\,,
    \label{e:lateorder_intro}
\end{equation}
where $\chi(z) = 0$ at the singular point in the leading-order term $V_0$ ($z = 0$) and $\gamma$ is a constant to be determined. By substituting the ansatz  (\ref{e:lateorder_intro})
into the recurrence relation \eqref{V.recur} and matching in the limit that $n \to \infty$, we find at the first two orders:
\begin{align}
\label{eq.Vchieq}\mathcal{O}(V_{n+1}):& \quad z = 2 \mu \diff{\chi}{z},\\
\label{eq.Vprefeq}\mathcal{O}(V_{n}):& \quad z \diff{G}{z} + G = \frac{\i}{z}\diff{\chi}{z} + 4 \mu \diff{\chi}{z}\diff{G}{z} + 2 \mu G \diff{^2\chi}{z^2}.
\end{align}
Using the condition $\chi(0) = 0$ and solving both of these equations gives
\begin{equation}
\chi = \frac{z^2}{4\mu}, \qquad G = \Lambda z^{-\i/2\mu}, 
\end{equation}
where $\Lambda$ is a constant.

\begin{figure}
\centering
\includegraphics[]{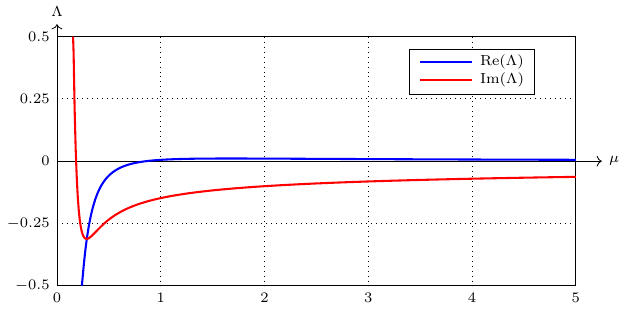}
\caption{The result of approximating $\Lambda$ different values of $\mu$. This approximation was obtained by evaluating the ratio in \eqref{eq.Vcomp} at each value of $\mu$ for $n = 1000$. This choice ensured that all values converged up to at least six decimal places. The value of $|\Lambda|$ becomes large as $\mu \to 0$.}\label{Fig.Lambda2}
\end{figure}

Determining the constant $\Lambda$ requires using matched asymptotic expansions to ensure that the late-order terms  (\ref{e:lateorder_intro}) are consistent with a local expansion of the solution in a small neighbourhood around the pole ($z = \mathcal{O}(\eps)$ as $\eps \to 0$), identified as the region in which the series terms of \eqref{eq.Vser} become comparable in size in the asymptotic limit, and the expansion therefore breaks down. From this matching, we can determine $\Lambda$, which depends on $\mu$. The technical details of this calculation are presented in Appendix \ref{s.lambda}. Computed values are presented in Figure \ref{Fig.Lambda2}. It is significant that $|\Lambda|$ appears to grow indefinitely as $\mu \to 0^+$; this suggests that our analysis breaks down in this limit, as the asymptotic scalings are different to those of the small-$\mu$ regime, studied in \cite{Chapman_shock_caustic}. The two analyses, ours and that provided in \cite{Chapman_shock_caustic}, combine to explain the solution behaviour for both $|\mu| \ll 1$ and $\mu = \mathcal{O}(1)$.

The late-order terms  (\ref{e:lateorder_intro}) therefore take the form
\begin{equation}\label{eq.VLOT2}
V_n \sim \frac{\Lambda z^{-\i/2\mu} \Gamma(n + \gamma)}{(z^2/4\mu)^{n+\gamma}} \quad \mathrm{as} \quad n \to \infty.
\end{equation}
For this expression to be consistent with the leading-order \eqref{eq.V0}, the strength of the singularity at $z=0$ in \eqref{eq.VLOT2} with $n=0$ must be the same as the strength of the singularity in $V_0$ (ie. equal to one). Applying this condition shows that $\gamma = \tfrac{1}{2}-\tfrac{\i}{4\mu}$. 

\subsection{Stokes switching}

Given the divergent asymptotic expansion (\ref{eq.Vser}) and the late-order terms (\ref{e:lateorder_intro}), it follows that an optimally truncated series has an exponential remainder of the form (\ref{e:stokes_intro}) (see (\ref{eq.VRNS}) below).  We can determine the location of Stokes curves, which are found where this exponential term is maximally subdominant compared to the leading-order term in the algebraic series \cite{Dingle}; in other words, where
\begin{equation}
    \mathrm{Im}(\chi)=0
    \quad\mathrm{and} 
    \quad\mathrm{Re}(\chi)>0.
\label{e:Stokes_def} 
\end{equation}
These are the curves across which the exponential term is switched by the algebraic series, and hence behaviour on an exponentially small scale can appear or vanish as these curves are crossed. The condition on $\mathrm{Re}(\chi)$ forces the changing contribution to be exponentially small in the asymptotic limit, rather than exponentially large.

Anti-Stokes curves are found where the exponential term becomes asymptotically comparable in size to the algebraic power series. These curves generally denote a change in dominance between two asymptotic contributions, and satisfy 
\begin{equation}
    \mathrm{Re}(\chi) = 0.
    \label{eq.antistokes}
\end{equation}
If an anti-Stokes curve is crossed, the terms in the asymptotic series reorder, as the exponential terms are comparable in size to the algebraic terms. The power series approximation \eqref{eq.Vser} is no longer valid in the region of the complex plane beyond the anti-Stokes curve. The Stokes structure is illustrated in Figure \ref{fig.zstokes}. 

Note that Figure \ref{fig.zstokes} is presented in terms of $\xi$. While the exponential asymptotic calculations were performed in terms of $z$, we recall that taking $\eps \to 0$ with $z = \mathcal{O}(1)$ corresponds to the limit $|\xi| \to \infty$. For simplicity of interpretation, we show the Stokes structure in terms of $\xi$. We will use the fact that $\mathrm{arg}(z) = \mathrm{arg}(\xi)$ to describe the Stokes structure in this section in terms of $\xi$.

In the red shaded regions of Figure \ref{fig.zstokes}, the exponential term is larger than the algebraic series. In the blue shaded regions, the exponential is smaller than the algebraic series. In the unshaded regions, there are no exponential contributions in the solution. 

In the region $\mathrm{arg}(\xi) \in (-3\pi/4,-\pi/4)$, which contains the negative imaginary axis, any possible exponential term would necessarily dominate the algebraic leading-order solution in the asymptotic limit, as $\mathrm{Re}(\chi) < 0$. We therefore conclude that there are no exponential contributions in this region, otherwise the boundary condition \eqref{eq.farfield} cannot be satisfied. This region is bounded by an anti-Stokes curve that satisfies \eqref{eq.antistokes} along $\mathrm{arg}(\xi) = -\pi/4$.

In the region $\mathrm{arg}(\xi) \in (-\pi/4, \pi/4)$, it can be seen that $\mathrm{Re}(\chi) > 0$, and  any exponential term present in the solution is small in the asymptotic limit. This region contains a Stokes curve that satisfies (\ref{e:Stokes_def}) along the positive real axis; the exponentially small behaviour in the solution must therefore appear as this curve is crossed into the upper half-plane.

Finally, the curves $\mathrm{arg}(\xi) = \pi/4$ and $3\pi/4$ are anti-Stokes curves that satisfy \eqref{eq.antistokes}, across which the exponential behaviour becomes comparable in size to the algebraic power series. The series approximation \eqref{eq.Vser} is no longer valid in this case, as the terms in the asymptotic series reorders as these curves are crossed. We will use transasymptotic analysis in Section \ref{sec.trans} to study the asymptotic solution beyond this anti-Stokes curve, and to approximate the location of the poles, which are generally located in the region $\mathrm{arg}(\xi) \in(\pi/4, 3\pi/4)$.

\subsection{Matched asymptotic expansions}
There are a number of different methods for determining the behaviour that is switched on as the Stokes curve is crossed, such as hyperasymptotics or Borel summation methods (eg. in \cite{BerryHowls1990,Berry1991,olde1995hyperasymptotic,olde1995hyperasymptotic2,olde2005hyperasymptotics,Sauzin-summability}). We apply the matched asymptotic expansion method of Ref.~\cite{Daalhuis}. This approach is typical of Stokes switching analysis (eg. see similar analyses in \cite{chapmanmortimer_2005,Chapman2006,Chapman2009,Chapman_shock_caustic,deng2023exponential,lustri2012,Lustri}) but we include an outline of the details here in order to showcase our full methodology.

We first truncate the power series \eqref{eq.Vser} after $N$ terms to minimise the truncation error. We apply the heuristic from \cite{Boyd1999} and truncate after the smallest term in the series. For this purpose, we use the late-order ansatz  (\ref{e:lateorder_intro}) to identify the value of $N$ at which consecutive terms of the series are the same size in the asymptotic limit that $\eps \to 0$ and $n \to \infty$.  This gives $N \sim |\chi|/\epsilon^2$ as $\epsilon \to 0$. We therefore set the truncation point to be $N = |\chi|/\epsilon^2 + \omega$, where $\omega \in [0,1)$ is chosen so that $N$ takes integer value. The truncated series is
\begin{equation}\label{eq.VserN}
V(z) = \sum_{n = 0}^{N-1} \epsilon^{2n}V_n(z) + R_N(z),
\end{equation}
where $R_N(z)$ is the truncation remainder. If the series is truncated optimally, this remainder will describe the exponentially small behaviour in the solution \cite{Berry1991}. We substitute the truncated series \eqref{eq.VserN} into the governing equation \eqref{eq.zeODE} and apply the recurrence relation \eqref{V.recur} to simplify. This produces  
\begin{equation}\label{eq.VRN}
\epsilon^2\left(2\mu \diff{R_N}{z} - 2 V_0 R_N + \ldots\right) + z R_N + \frac{\i}{2} \sim  z \eps^{2N} V_N\quad \as \quad \epsilon \to 0,
\end{equation}
where the omitted terms are subdominant to those which were retained in the limit $\eps \to 0$.

The right-hand side of (\ref{eq.VRN}) is exponentially small compared to the left-hand side except in an asymptotically small neighbourhood surrounding the Stokes curve. Applying a Liouville-Green (or WKB) ansatz to the homogeneous version of \eqref{eq.VRN} shows that, away from the Stokes curve
\begin{equation}
R_N \sim k G \e^{-\chi/\epsilon^2} \quad \as \quad \eps \to 0,
\end{equation}
where $k$ is some arbitrary constant. Motivated by this expression, we apply a variation of parameters approach to describe the solution to the inhomogeneous expression \eqref{eq.VRN}. We replace to constant $k$ with a function $\mathcal{S}(z)$, known as a Stokes multiplier, to give
\begin{equation}\label{eq.VRNS}
R_N = \mathcal{S}G \e^{-\chi/\epsilon^2}.
\end{equation}
The Stokes multiplier $\mathcal{S}$ captures the rapid variation in the exponential behaviour in the vicinity of the Stokes curve, while being asymtotically constant away from the Stokes curve. Substituting this expression into \eqref{eq.VRN} and using the late-order ansatz  (\ref{e:lateorder_intro}) for the right-hand side gives
\begin{equation}
\epsilon^2 \diff{\mathcal{S}}{z} G\e^{-\chi/\epsilon^2} \sim \frac{z \eps^{2N} }{2\mu}\frac{ G\Gamma(N + \gamma)}{\chi^{N+\gamma}}\quad \as \quad \eps \to 0.
\end{equation}
We apply a change of variables to use $\chi$ as the independent variable, and use \eqref{eq.Vchieq} to simplify the expression. After some algebra, this gives
\begin{equation}
\diff{\mathcal{S}}{\chi} \sim \eps^{2N-2}\frac{\Gamma(N+\gamma)}{\chi^{N+\gamma}}\e^{\chi/\epsilon^2} \quad \as \quad \eps \to 0.
\end{equation}
We now express the singulant in polar coordinates, $\chi = r \e^{\i\theta}$. The direction of rapid variation is perpendicular to the Stokes curve \cite{Daalhuis}, which corresponds to setting $r$ to be constant and considering only angular variation. We therefore set $N = r/\epsilon^2 + \omega$. This gives
\begin{equation}
\diff{\mathcal{S}}{\theta} \sim \i r\e^{\i\theta}\eps^{2r/\epsilon^2 +2\omega - 2}\frac{\Gamma(r/\epsilon^2 + \omega + \gamma)}{\chi^{r/\epsilon^2 + \omega+ \gamma}}\exp\left(\frac{r}{\eps^2}\e^{\i\theta}\right)\quad \as \quad \eps \to 0.
\end{equation}
Applying Stirling's formula and simplifying this expression gives
\begin{equation}
\diff{\mathcal{S}}{\theta} \sim \frac{\i\sqrt{2\pi r}}{\eps^{2\gamma+1}}\exp\left(\frac{r}{\eps^2}\left(\e^{\i\theta}-1-\i\theta\right) + \i\theta(1-\omega-\gamma)\right)\quad \as \quad \eps \to 0.
\end{equation}
The right-hand side of this expression is exponentially small except in a small neighbourhood near $\theta = 0$, which corresponds to the Stokes curve. To determine the variation near this curve, we define a local variable $\theta = \eps \vartheta$, and find that
\begin{equation}
\diff{\mathcal{S}}{\vartheta} \sim \frac{\i\sqrt{2\pi r}}{\eps^{2\gamma}}\exp\left(-\frac{r\theta^2}{2}\right)\quad \as \quad \eps \to 0.
\end{equation}
Integrating this expression gives
\begin{equation}
\mathcal{S} \sim \int_{-\infty}^{\vartheta/\sqrt{r}}\e^{-s^2/2}\d s + C,
\end{equation}
where $C$ is some constant of integration. In fact, taking the limit that $\vartheta \to -\infty$ corresponds to moving into the sector $\arg(\xi) < 0$, in which the exponential is not present. Hence $C = 0$. This gives the jump in exponential behaviour as the Stokes curve is crossed as
\begin{equation}
\left[\mathcal{S}\right]_-^+ = \lim_{{\vartheta} \to \infty} \mathcal{S} - \lim_{{\vartheta} \to -\infty} \mathcal{S} \sim \frac{2\pi \i}{\eps^{2\gamma}},
\end{equation}
and therefore the exponential contribution for $\arg(\xi) > 0$ is given by
\begin{equation}
V_{\mathrm{exp}} \sim \frac{2\pi \i G }{\eps^{2\gamma}} \e^{-\chi/\eps^2}\quad \as \quad \eps \to 0.
\end{equation}
Returning to original coordinates using $\xi = z/\eps$ and $U = \eps V$ gives
\begin{equation}\label{eq.Uexp}
U_{\mathrm{exp}} \sim 2\pi\i\Lambda \xi^{-\i/2\mu}\e^{-\xi^2/4\mu}\quad \as \quad |\xi| \to \infty.
\end{equation}
This expression is the leading-order behaviour behaviour of the exponential terms in regions where these terms are present but small in the asymptotic limit (ie. Sectors \ding{173} and \ding{175} in Figure \ref{fig.zstokes}). Note that the $\eps$ parameter disappears entirely from this expression. This is expected, as it was an artificial parameter introduced to index the size of $|\xi|$. 

We have therefore determined the solution behaviour in the complex plane in Sectors \ding{173} and \ding{175} from Figure \ref{fig.zstokes}, or $\mathrm{arg}(\xi) \in (0,\pi/2)$ or $(3\pi/4,\pi)$. To continue beyond from these regions into Sector \ding{174}, or $\mathrm{arg}(\xi) \in (\pi/4,3\pi/4)$ requires understanding how the exponential terms behave as the anti-Stokes curve is crossed.

\section{Transasymptotic Analysis}\label{sec.trans}

We wish to determine the location of the poles and zeros in the solution to \eqref{eq.z0ODE}--\eqref{eq.farfield}, which requires studying the behaviour of the solution as the anti-Stokes curves on $\mathrm{arg}(\xi) = \pi/4$ and $3\pi/4$ are crossed. As these curves are crossed, the exponentially-small terms that appear across the Stokes curves, described in Section \ref{sec.Stokes}, grow and become asymptotically large. This changes the asymptotic balance of the solution, and the series approximation in \eqref{eq.Vser} is no longer valid.  Poles are typically located along curves that asymptotically approach anti-Stokes curves in the large-$|\xi|$ limit \cite{costin2002formation}, so locating the poles requires understanding the behavior of the solution in the region where the exponentials grow and become asymptotically large. 

Approximating the solution behaviour in this region requires first writing the solution in the region where the exponential terms are small as a transseries, or a double series in multiple scales which describes both algebraic and exponential asymptotic contributions to the solution behaviour. The final step is to sum the transseries over all exponential terms at each algebraic order so that it remains valid even where the exponential terms reorder in size. This process is known as transasymptotic analysis, and it will produce an expression which we study to determine the location of poles caused by growing exponential terms in the asymptotic behaviour.

\subsection{Formulating the transseries}\label{sec:formulatetrans}

We begin by integrating \eqref{eq.z0ODE} and applying the boundary condition in \eqref{eq.farfield}, to obtain 
\begin{equation}
    \label{eq.z2ODE}
    \mu \diff{U}{\xi} =- \frac{1}{2}\xi U + \frac{1}{2}U^2 - \frac{\i}{4}.
\end{equation}
This is not strictly necessary, but it reduces the order of the differential equation and makes the subsequent analysis more straightforward.

The transseries must contain an algebraic power series for $U$, as well as the exponentially small terms whose leading-order behaviour is given by \eqref{eq.Uexp}. Furthermore, due to the nonlinearity in the problem, the transseries will contain exponential terms with exponents that are integer multiples of the exponential term in \eqref{eq.Uexp}  \cite{Sauzin-summability,ANICETO20191}. The presence of these additional exponentials suggests that we apply a transseries of the form 
\begin{equation}\label{eq.transseries}
U(\xi) \sim \sum_{n=0}^{\infty} \sigma^n \e^{-n \xi^2/{4\mu}} U_n(\xi) \qquad \mathrm{as} \qquad |\xi| \to \infty,
\end{equation}
where $U_n(\xi)$ can be represented as an algebraic series in $\xi$, and $\sigma$ is a transseries parameter that will be determined from the exponential contribution derived in Section \ref{sec.Stokes}, given explicitly in \eqref{eq.Uexp}. 

In Sector \ding{172} of Figure \ref{fig.zstokes}, corresponding to $\mathrm{arg}(\xi) \in(-\pi, 0)$, there are no exponentially small contributions. We capture this behaviour by setting $\sigma = 0$, which means that the asymptotic solution behaviour is entirely described by the power series \eqref{eq.F0ser}. The analysis in Section \ref{sec.Stokes} showed that exponentially small terms are present in sectors \ding{173} and \ding{175}, or $\mathrm{arg}(\xi) \in (0, \pi/4)$ and $(3\pi/4,\pi)$ respectively, which have the leading-order behaviour given in \eqref{eq.Uexp}. In these sectors we set $\sigma \neq 0$, which produces exponentially small transseries terms.

The power series expansion for $U_0$ as $\xi \to \infty$ may be found by direct substitution (or by matching with the outer expansion via (\ref{eq:3termfarfield})); we find that 
\begin{equation}\label{eq.F0ser}
U_0(\xi) \sim -\frac{\i}{2\xi} + \left(-\frac{1}{4} - \i\mu\right)\frac{1}{\xi^3} +  \left(\frac{\i}{4} - \frac{5\mu}{2} - 6 \i \mu^2 \right)\frac{1}{\xi^5} + \mathcal{O}(\xi^{-7})\quad \mathrm{as} \quad |\xi| \to \infty.
\end{equation} 
Matching at higher orders of the exponential term gives equations for the subsequent transseries components,
\begin{equation}\label{eq.Uneq}
    \mu \diff{U_n}{\xi} = \frac{(n-1)}{2}\xi U_n+ \frac{1}{2}\sum_{j=0}^{n}U_j U_{n-j},\qquad n \geq 1.
\end{equation}
We treat \eqref{eq.Uneq} differently for the case $n = 1$, which causes the first term on the right-hand side to vanish. A Green-Liouville (or WKB) analysis of the solution for $n = 1$ gives
\begin{equation}
    U_1 \sim K_1 \xi^{-\i/2\mu} \qquad \mathrm{as} \qquad |\xi| \to \infty.
\end{equation}
where $K_1$ is a constant that has yet to be determined. If $n > 1$, the first term on the right-hand side of \eqref{eq.Uneq} does not vanish, and the asymptotic dominant balance is obtained by balancing the terms on the right-hand side. This gives
\begin{equation}\label{eq.WKBresult}
    U_n \sim K_n \xi^{1-n(1+\i/2\mu)} \qquad \mathrm{as} \qquad |\xi| \to \infty,
\end{equation}
where $K_n = (-1)^{n+1}K_1$. The power of $\xi$ in \eqref{eq.WKBresult} suggests that we write the transseries \eqref{eq.transseries} as
\begin{equation}\label{eq.Fseries1}
U(\xi) \sim \sum_{m=1}^{\infty}\frac{a^{(0)}_m}{\xi^{2m-1}} + \sum_{n=1}^{\infty} \sigma^n \xi^{-n(1 + \mathrm{i}/2\mu)}  \e^{-n \xi^2/4\mu} \sum_{m=0}^{\infty} \frac{a^{(n)}_m}{\xi^{2m-1}},
\end{equation}
where $a^{(n)}_m$ is the constant in the $m$th term in the power series associated with $U_n$. The range of index $m$ is chosen for $a_m^{(0)}$ so that the powers of $\xi$ are consistent for corresponding choices of $m$ across the $n = 0$ and $n > 0$ cases. Finally, we are able to write \eqref{eq.Fseries1} as a typical transseries ansatz by defining a new variable $\tau$, 
\begin{equation}\label{eq.tau}
\tau = \sigma  {\xi}^{-(1+\i/2\mu)} \e^{-\xi^2/4\mu},
\end{equation}
Note that $\tau$ is exponentially small in sectors \ding{173} and \ding{175} of the solution, corresponding to regions in which the exponential terms are present and small. We collect the terms to write 
\begin{equation}\label{eq.transseries5}
U(\xi) \sim  \sum_{n=0}^{\infty} \tau^n \sum_{m=0}^{\infty} \frac{a^{(n)}_m}{\xi^{2m-1}}.
\end{equation}
which is valid in the limit that $|\xi| \to \infty$ in regions where the exponential contributions are asymptotically small, or $\tau \to 0$; this corresponds to $|\xi| \to \infty$ in the sector $\mathrm{arg}(\xi) \in (-5\pi/4,\pi/4)$. 

We require that $a^{(0)}_0 = 0$ to be consistent with \eqref{eq.Fseries1} and the algebraic series solution, which is $\mathcal{O}(\xi^{-1})$ as $|\xi| \to \infty$. The values of $a^{(0)}_n$ for $n \geq 1$ may be read directly from the series \eqref{eq.F0ser}, giving
\begin{equation}\label{eq.a0m}
a^{(0)}_1 = -\frac{\i}{2}, \qquad a^{(0)}_2 = -\frac{1}{4}-\i\mu, \qquad a^{(0)}_3 = \frac{\i}{4}-\frac{5\mu}{2} - 6 \i\mu^2.
\end{equation}
We may obtain all terms of the form $a^{(0)}_m$ by directly calculating subsequent terms of the algebraic power series \eqref{eq.F0ser}. 

We determine $\sigma$ by considering the term indexed by $a^{(1)}_0$, given by $\tau \xi a^{(1)}_0$. This term is the leading-order behaviour of the first exponential correction to the algebraic series, and is given by $U_{\mathrm{exp}}$ from \eqref{eq.Uexp} in sectors \ding{173} and \ding{175}. By equating these expressions, we have
\begin{equation}\label{eq.tauint}
\tau \xi a^{(1)}_0 = 2\pi \i \Lambda \xi^{-\i/2\mu}\e^{-\xi^2/4\mu}.
\end{equation}
Substituting the definition of $\tau$ in \eqref{eq.tau} into \eqref{eq.tauint} and simplifying the result shows
\begin{equation}
    \sigma a^{(1)}_0 = 2\pi\i \Lambda.
\end{equation}
We are free to set $\sigma = 2\pi\i\Lambda$, which gives $a^{(1)}_0 = 1$. Recall that $\sigma = 0$ in sector \ding{172}. The value of $\sigma$ is different on either side of Stokes curves in the complex plane -- this is typical behaviour for the transseries parameter.

We now substitute the transseries \eqref{eq.transseries5} into the differential equation \eqref{eq.z2ODE} and match at higher orders in $\tau$, recalling that $\tau \to 0$ exponentially fast as $|\xi| \to \infty$. At higher orders, we find 
\begin{equation}\label{eq.a11a21}
\mathcal{O}(\tau): \, a_1^{(1)} = \frac{\i}{2} + \frac{1}{8\mu},\qquad 
\mathcal{O}(\tau \xi^{-1}): \, a_2^{(1)} = \frac{3\i \mu}{2}  + \frac{1}{2} + \frac{1}{128\mu^2}.
\end{equation}
We can continue to match in this fashion, and it is straightforward to compute an arbitrary number of coefficients $a_m^{(n)}$ in both the $m$ and $n$-indexed directions. Incrementing $m$ gives
\begin{equation}\label{eq.a13}
 a_3^{(1)} = 10\i \mu^2+\frac{55\mu}{12} -\frac{7\i}{16}+\frac{1}{24\mu}-\frac{\i}{256\mu^ 2}  + \frac{1}{3072\mu^3},
\end{equation}
while incrementing $n$ gives
\begin{equation}\label{eq.a02}
 a_0^{(2)} = -1, \qquad a_1^{(2)} =2\mu-\frac{1}{4\mu},\qquad a_2^{(2)}=-12\mu^2-7\i \mu+\frac{1}{4}+\frac{\i}{8\mu}-\frac{1}{32\mu^2}.
\end{equation}
The values presented in \eqref{eq.a0m} and \eqref{eq.a11a21} give sufficient information about the behaviour in sectors \ding{173} and \ding{175} to approximate the pole locations in sector \ding{174}. 

\subsection{Transasymptotic analysis}\label{sec:transanalysis}

The transseries \eqref{eq.transseries5} is a valid asymptotic expression for the solution behaviour in $\mathrm{arg}(\xi) \in (-5\pi/4,\pi/4)$. In sector \ding{172}, or $\mathrm{arg}(\xi) \in (-\pi, \pi)$, the exponential terms are absent, so $\sigma = 0$. In sectors \ding{173} and \ding{175}, or $\mathrm{arg}(\xi) \in (0, \pi/4)$ and $\mathrm{arg}(\xi) \in (-5\pi/4, -\pi)$, the exponential terms are small, and hence $\tau$ is exponentially small in the limit that $|\xi| \to \infty$. To determine the location of poles in the solution, we must approximate the behaviour in $\mathrm{arg}(\xi) \in (\pi/4, 3\pi/4)$, or sector \ding{174}. Continuing the solution into this region requires extending the expression \eqref{eq.transseries5} so that it remains valid as the anti-Stokes curves are crossed and the exponential terms in the transseries become asymptotically large.

The key idea of transasymptotic analysis is that the terms in the transseries \eqref{eq.transseries} can be resummed to approximate the behaviour in regions where $\tau$ is not exponentially subdominant compared to $\xi$. This is possible because the transseries \eqref{eq.transseries} describes all of the subdominant exponential terms in the solution; hence, any change in the asymptotic balance of the problem corresponds to a rearrangement of the terms in the transseries.

For convenience, we rewrite the transseries as 
\begin{equation}\label{eq.transseries2}
U(\xi) \sim  \sum_{n=0}^{\infty} \sum_{m=0}^{\infty} \frac{\tau^n a^{(n)}_m}{\xi^{2m-1}},
\end{equation}
as $|\xi| \to \infty$ and $\tau \to 0$, where we set $a_0^{(0)} = 0$ to be consistent with \eqref{eq.transseries5}. This is a divergent asymptotic series, and we change the order of summation to obtain the asymptotically valid expression
\begin{equation}\label{eq.transseries3}
U(\xi) \sim  \sum_{m=0}^{\infty} \frac{1}{\xi^{2m-1}}\left(\sum_{n=0}^{\infty}\tau^n a^{(n)}_m\right)
\end{equation}
as $|\xi| \to \infty$ and $\tau \to 0$. Finally, we sum the series over the index $n$. This is the resummation step, which allows us to describe the asymptotic solution behaviour even in regions where $\tau$ is not small. We write the resultant expression as
\begin{equation}\label{eq.transseries4}
U(\xi) \sim  \sum_{m=0}^{\infty} \frac{A_m(\tau)}{{\xi^{2m-1}}}, 
\end{equation}
which is valid as $|\xi|\to\infty$, irrespective of the size of $\tau$. In regions of the complex $\xi$-plane where the inner sum of \eqref{eq.transseries3} converges, it is equal to $A_m(\tau)$; however, $A_m(\tau)$ can also be evaluated in regions of the complex plane where the inner sum does not converge. It is a continuation of the inner sum from \eqref{eq.transseries3} to the entire complex plane. 

Substituting \eqref{eq.transseries4} into the differential equation \eqref{eq.z2ODE} gives an expression for the series terms
\begin{equation}\label{eq.transsub}
    \mu \sum_{m=0}^{\infty} \left(\frac{1}{\xi^{2m-1}}\diff{A_m}{\xi}  - \frac{(2m-1) A_m}{\xi^{2m}}\right) = -\frac{\xi}{2}\sum_{m=0}^{\infty} \frac{A_m}{{\xi^{2m-1}}} + \frac{1}{2} \sum_{l=0}^{\infty}\sum_{m=0}^{\infty} \frac{A_l A_m}{{\xi^{2l+2m-2}}}-\frac{\mathrm{i}}{4}.
\end{equation}
Because $\tau$ is a function of $\xi$, the chain rule gives
\begin{equation}\label{eq.dAmdx}
\diff{A_m}{\xi} = \diff{\tau}{\xi}\diff{A_m}{\tau} 
= -\tau \left[\frac{\xi}{2\mu} + \frac{1}{\xi}\left(1 + \frac{\i}{2\mu}\right)\right] \diff{A_m}{\tau}. 
\end{equation}
Applying \eqref{eq.dAmdx} to \eqref{eq.transsub} and matching terms as $|\xi| \to \infty$ gives the recurrence relation 
\begin{align}
\label{eq.A0ode}A_0 - A_0^2 - \tau \diff{A_0}{\tau} &= 0,\\
\label{eq.A1ode}A_1 + 2\mu A_0 -2 A_0 A_1  - \tau \diff{A_1}{\tau} - \left(\i  + 2\mu \right) \tau \diff{A_0}{\tau} + \frac{\i}{2} & = 0,\\
\label{eq.Akode}A_{k+1} -2\mu \left(2k-1\right) A_k  - \sum_{m=0}^{k+1} A_m A_{k+1-m}  - \tau \diff{A_k}{\tau} - \left(\i  + 2\mu \right) \tau \diff{A_k}{\tau}  & = 0, \qquad k \geq 1.
\end{align}
Solving \eqref{eq.A0ode} gives
\begin{equation}\label{eq.A0k}
A_0(\tau) = \frac{\tau}{k + \tau },
\end{equation}
where $k$ is an arbitrary constant. To determine the value of this constant, we must match this expression to the leading exponential behaviour in sector \ding{173}, or \eqref{eq.Uexp}, where $\tau$ is small. In this region, the inner sum from \eqref{eq.transseries3} converges, and 
\begin{equation}\label{e.A0trans}
A_0(\tau) = \sum_{n=0}^{\infty} \tau^n a_0^{(n)} \sim a_0^{(0)} +\tau  a_0^{(1)} + \mathcal{O}(\tau^2) \quad \mathrm{as} \quad \tau \to 0.
\end{equation}
Recall that $a_0^{(0)} = 0$ and $a_0^{(1)} = 1$. Hence, we find that $A_0(\tau) \sim \tau$ as $\tau \to 0$, so $k=1$ and
\begin{equation}\label{eq.A0}
A_0(\tau) = \frac{\tau}{1 +  \tau}.
\end{equation}
This rational expression is singular at $\tau = -1$. This is sufficient for us to approximate the pole locations to leading order; however, it is straightforward to improve the approximation by calculating the first correction term, $A_1$.

The governing equation for $A_1(\tau)$ is given by \eqref{eq.A1ode}, and the boundary condition can be found by considering the behaviour of $A_1(\tau)$ in sector \ding{173}, where $\tau \to 0$ and the inner sum of \eqref{eq.transseries3} converges. In this region, we have
\begin{equation}
A_1(\tau) = \sum_{n=0}^{\infty} \tau^n a_1^{(n)} \sim a_1^{(0)} +  \tau  a_1^{(1)}+ \mathcal{O}(\tau^2) \quad \mathrm{as} \quad \tau \to 0.
\end{equation}
Using the values determined in \eqref{eq.a11a21}, we find that
\begin{equation}
A_1(\tau) \sim 1 + \left( \frac{\i}{2} + \frac{1}{8\mu}\right)\tau \quad \mathrm{as} \quad \tau \to 0.
\end{equation}
Solving \eqref{eq.A1ode} subject to this boundary condition gives
\begin{equation}\label{eq.A1}
A_1(\tau) = \frac{1}{(1+\tau)^2}\left[\frac{\tau(1 - 4\mu \i)(1 + 4\mu \i \tau)}{8\mu} -\frac{\i}{2}\right].
\end{equation}
The functions $A_0(\tau)$ and $A_1(\tau)$ provide sufficient information for us to approximate the pole locations in sector \ding{174} up to the first correction term in the limit that $|\tau| \to 0$. We could continue this process to find further corrections using \eqref{eq.Akode}, with boundary conditions that depend on $a_m^{(0)}$ and $a_m^{(1)}$. For instance, the next correction term is
\begin{align}\nonumber
A_2(\tau)=&-\frac{1}{128 \mu^2(1+\tau^3)}\left(32 \mu ^2 \left(32 \mu ^2+12 i \mu -1\right) \tau ^3+32 \mu ^2 (1+4 i \mu )+\right.\\ &+\left.\left(192 i \mu ^3+32 \mu ^2-1\right) \tau +\left(1536 \mu ^4+704 i \mu ^3-128 \mu ^2-16 i \mu +1\right) \tau ^2\right).\label{e.A2}
\end{align}

Finally, we approximate the behaviour of $U(\xi)$ in sector \ding{174} using the resummed series \eqref{eq.transseries4}, as well as \eqref{eq.A0} and \eqref{eq.A1}, 
\begin{equation}\label{eq.Fs}
U(\xi) \sim \frac{\tau \xi}{1 + \tau} +   \frac{1}{\xi}\frac{1}{(1+\tau)^2}\left[ \frac{\tau(1 - 4\mu \i)(1 + 4\mu \i \tau)}{8\mu} - \frac{\i}{2}\right] + \mathcal{O}(|\xi|^{-3}) \quad \mathrm{as} \quad |\xi| \to \infty.
\end{equation}
This expression has the correct asymptotic behaviour in sector \ding{173}, and can be extended past the anti-Stokes curves at $\mathrm{arg}(\xi) = \pi/4$ and $3\pi/4$.

Note that (\ref{eq.Fs})is only valid in regions where $\sigma \neq 0$, or the upper-half complex plane. In the lower-half plane we instead have $\sigma = 0$, and therefore $A_k(\tau) = 0$ for $k \geq 0$. We therefore do not predict any far-field poles to be located in sector \ding{172}, or the entire lower-half complex plane (poles close to the origin of the $\xi$-plane may lie in the lower-half plane, and indeed this does occur for $\mu\lessapprox 0.147$ \cite{vandenheuvel2022burgers}).

\section{Pole Locations}\label{sec:polelocations}

From the form of \eqref{eq.A0}, we see that that the singularities occur at $\tau \sim -1$ as $|\xi| \to \infty$, which can be mapped back to coordinates in $(x,t)$. We can obtain a correction term to this expression which permits us to approximate the pole locations more accurately. We will write the corrected pole location as 
\begin{equation}
\tau \sim \tau_0 + \tau_1 \xi^{-2} + \tau_2 \xi^{-4} + \ldots \quad \mathrm{as} \quad |\xi| \to \infty.
\end{equation}
We know that all poles in the solution must be simple, and therefore write 
\begin{equation}\label{eq.Fs0}
\frac{U(\xi)}{\xi} = \frac{P(\tau)}{\tau - \tau_0 - \tau_1 \xi^{-2} - \tau_2\xi^{-4} - \ldots},
\end{equation}
where $P(\tau)$ is a polynomial expression. Expanding \eqref{eq.Fs0} as $|\xi| \to \infty$ gives
\begin{equation}\label{eq.Fs1}
\frac{U(\xi)}{\xi} \sim \frac{P(\tau)}{\tau - \tau_0} + \frac{\tau_1 P(\tau)}{\xi^2(\tau - \tau_0)^2} + \mathcal{O}(|\xi|^{-4}).
\end{equation}
Now, we compare \eqref{eq.Fs1} with \eqref{eq.Fs}. Dividing \eqref{eq.Fs} throughout by $\xi$ gives
\begin{equation}\label{eq.Fss}
\frac{U(\xi)}{\xi} \sim\frac{\tau }{1+\tau} +   \frac{1}{\xi^2 (1+\tau)^2}\left[\frac{\tau(1 - 4\mu \i)(1 + 4\mu \i \tau)}{8\mu} -\frac{\i}{2}\right]   \quad \mathrm{as} \quad |\xi| \to \infty.
\end{equation}
Matching the first term of \eqref{eq.Fs1} with that of \eqref{eq.Fss} near $\tau = \tau_0$ gives
\begin{equation}
\tau_0 = -1,\qquad P(\tau) = \tau.
\end{equation}
We obtain $\tau_1$ by matching the second term of \eqref{eq.Fs1} with that of \eqref{eq.Fss} in the vicinity of $\tau = \tau_0$ for large $|\xi|$. Matching \eqref{eq.Fs1} and \eqref{eq.Fss} requires 
\begin{equation} 
P(\tau_0)\tau_1=\left[\frac{\tau_0(1 - 4\mu \i)(1 + 4\mu \i \tau_0)}{8\mu} -\frac{\i}{2}\right],
\end{equation}
and gives
\begin{equation}\label{eq.tau1}
\tau_1 = -\frac{\i}{2}+ \frac{1}{8\mu}-2\mu.
\end{equation}
To obtain further corrections to the pole position, we would have to include further corrections to $U(\xi)$ in \eqref{e.A0trans}. For instance, if we include $A_2(\tau)$ from \eqref{e.A0trans}, a similar analysis shows that
\begin{equation}\label{e.tau2}
    \tau_2=8 \mu^2-\frac{1}{128 \mu^2}+\frac{9 \mu \i }{2}.
\end{equation}

We now have sufficient information for us to approximate the pole locations explicitly in terms of $\xi$. We need to invert the expression
\begin{equation}
\tau = \sigma \xi^{-\alpha}\e^{-\xi^2/4\mu} \sim -1 + \tau_1 \xi^{-2} + \ldots \quad \mathrm{as} \quad |\xi| \to \infty,
\end{equation}
where we define $\alpha = 1 + \i/2\mu$ for notational convenience, and $\tau_1$ is given in \eqref{eq.tau1}. After taking logarithms of both sides and some further algebraic manipulation, we find
\begin{align}\label{e.logrel}
\log(\sigma) - \frac{\alpha}{2}\log(\xi^2) - \frac{\xi^2}{4\mu} + (2  M + 1)\pi \i \sim \log\left(1 - \frac{\tau_1}{\xi^2} + \ldots\right)  ,
\end{align}
where $M$ is any integer. For large values of $|M|$, the poles will be located at 
\begin{equation}\label{eq.poleleading}
    \xi \sim \pm\sqrt{8 \mu \i M \pi} \qquad \mathrm{as} \qquad |M| \to \infty.
\end{equation}
Here the positive (negative) sign is associated with positive (negative) $M$, so that the poles lie in the first (second) quadrant.  The result is a string of poles along the rays $\mathrm{arg}(\xi) = \pi/4$ and $\mathrm{arg}(\xi) = 3\pi/4$, corresponding to the anti-Stokes curves labelled in Figure \ref{fig.zstokes}.  The expression (\ref{eq.poleleading}) agrees with the leading-order approximation from \cite{vandenheuvel2022burgers}, obtained from a transcendental equation derived using properties of the exact solution (\ref{eq:exactsoln}).  

The derivation of \eqref{eq.poleleading} only requires $A_0(\tau)$, which in turn requires only the explicit calculation of $a_0^{(0)}$ and $a_0^{(1)}$.  In order to improve the accuracy of (\ref{eq.poleleading}), we will use the correction term $A_1(\tau)$. We set
\begin{equation}
    r_\mathrm{p} = \log(\sigma) + (2M+1)\pi\i,\label{eq:r-poles}
\end{equation}
and we will find an asymptotic expression for the pole location in the limit that $|r_{\mathrm{p}}| \to \infty$, which corresponds to $M \to \infty$. Although we are only determining the asymptotic behaviour of the pole locations in this limit, we will compare our results to numerical solutions in Section \ref{sec.results} and find that the approximations are accurate even for poles that are relatively close to the origin.

Note that \eqref{e.logrel} only contains $\xi^2$ terms, so we can write a series expression for this variable. Due to the logarithmic term in \eqref{e.logrel}, we expect that this series also contains logarithmic terms. We set
\begin{equation}
    \xi^2 \sim c_0 r_{\mathrm{p}} + c_1 \log(4 \mu r_\mathrm{p}) + c_2 +  \frac{c_3 \log(4 \mu r_\mathrm{p})}{r_\mathrm{p}} + \frac{c_4}{r_\mathrm{p}} + \ldots,
\end{equation}
where the omitted terms are small in the limit that $|r_\mathrm{p}| \to \infty$, and the $4\mu$ in the logarithmic terms is included for algebraic convenience. By substituting this into \eqref{e.logrel} and matching as $|r_\mathrm{p}| \to \infty$, we find
\begin{equation}
    c_0 = 4\mu , \qquad c_1 = -2\mu\alpha, \qquad c_2 = 0,  \qquad c_3 = \alpha^2\mu, \qquad c_4 = \tau_1 .
\end{equation}
We combine these terms to obtain
\begin{equation}\label{e:poleloc}
    \xi^2 \sim 4\mu r_{\mathrm{p}}  -2\mu\alpha \log(4\mu r_{\mathrm{p}}) + \frac{\alpha^2 \mu  \log(4 \mu  r_{\mathrm{p}}) + \tau_1}{r_{\mathrm{p}}} +  \ldots,
\end{equation}
where $r_{\mathrm{p}}$ depends on $M$, which takes integer values that index the poles. The expression \eqref{e:poleloc} gives an asymptotic expression for $\xi^2$, and hence it is straightforward to compute $\xi$. More terms could be obtained using the expression in \eqref{e.tau2} for $\tau_2$, but we do not require any further corrections here. Note that this approximation becomes more accurate as $|M|$ grows, because $|r_{\mathrm{p}}|$ grows correspondingly.

\section{Zero Locations}\label{sec:zerolocations}

We apply an essentially identical method to determine the asymptotic location of zeroes of $U(\xi)$. We write the location of the zeroes of $U(\xi)$ as $\tau \sim \tau_0 + \tau_1 \xi^{-2} + \tau_2 \xi^{-4} + \ldots$ as $|\xi| \to 0$, as before. From \eqref{eq.A0}, it is evident that $\tau_0 = 0$ yields $A_0(\tau_0) = 0$. 

Substituting the series expression for $\tau$ into \eqref{eq.Fss} and setting $U(\xi) = 0$ gives 
\begin{equation}\label{eq.tauzer}
     \frac{\tau_1\xi^{-2} }{1+\tau_1\xi^{-2}} +   \frac{1}{\xi^2 (1+\tau_1\xi^{-2})^2}\left[\frac{\tau_1\xi^{-2}(1 - 4\mu \i)(1 + 4\mu \i \tau_1\xi^{-2})}{8\mu} -\frac{\i}{2}\right] + \mathcal{O}(|\xi|^{-4}) = 0,
\end{equation}
as $|\xi| \to \infty$. Matching terms in \eqref{eq.tauzer} at leading order as $|\xi| \to \infty$ gives $\tau_1 = \i/2$, and hence that the location of the zeroes of $U(\xi)$ satisfies $\tau \sim \i/2\xi^2$ as $|\xi| \to \infty$. We determine the location of zeroes by inverting the expression 
\begin{equation}\label{e.tauinv}
\tau = \sigma \xi^{-\alpha}\e^{-\xi^2/4\mu} \sim  \tau_1\xi^{-2}  + \ldots \quad \mathrm{as} \quad |\xi| \to \infty.
\end{equation}
By taking the logarithm of both sides and rearranging, we find that
\begin{align}\label{e.logrelbb}
\log(\sigma) - \frac{\alpha-2}{2}\log(\xi^2) - \frac{\xi^2}{4\mu} + 2  M \pi \i \sim \log\left(\frac{\mathrm{i}}{2} + \ldots\right), 
\end{align}
where $M$ is any integer. This approximation is valid for large $|\xi|$, and hence $M$. For large values of $|M|$, this expression gives the asymptotic relation
\begin{equation}\label{eq.zeroleading}
    \xi \sim \pm\sqrt{8 \mu \i M \pi} \qquad \mathrm{as} \qquad |M| \to \infty.
\end{equation}
which is identical to the leading behaviour of the pole locations, given in \eqref{eq.poleleading}. The description of the pole and zero locations in \eqref{eq.poleleading} and \eqref{eq.zeroleading} is consistent with the behaviour seen in \cite{vandenheuvel2022burgers}, where it was shown that the poles and zeroes interlace. Consequently the difference between the pole zero locations is small compared to $|\xi|$, and should not appear at leading order as $|M| \to \infty$. 

We use identical methods to the previous section to invert \eqref{e.logrelbb}. First, comparing this equation with \eqref{e.logrel} we see that it is natural to define 
\begin{equation}
    r_{\mathrm{z}}=\log(2\sigma) + \left(2M-\tfrac{1}{2}\right)\pi\i,\label{eq:r-zeros}
\end{equation}
and we then invert expression \eqref{e.logrel} to obtain
\begin{equation}\label{e:zeroloc}
        \xi^2 \sim 4\mu r_{\mathrm{z}} - 2(\alpha - 2)\mu \log(4\mu r_{\mathrm{z}}) - 6 \mu \pi \i +4\mu\log 2 \quad \mathrm{as} \quad |\xi| \to \infty.
\end{equation}
Note that $\pi\mathrm{i}/2-\log 2=\log(\tau_1)$. Because $\tau_1$ is the leading contribution to the position of the zeros, we do not have sufficient accuracy to compute $r_{\mathrm{z}}^{-1}$ corrections for the variable $\xi$, unlike in the pole calculations. We could obtain more terms of accuracy in this expression using the value of $\tau_2$ from \eqref{e.tau2} in the inverted expression for $\tau$ \eqref{e.tauinv}, but \eqref{e:zeroloc} is sufficient to obtain reasonable approximations for the zero locations. Finally, note that there is some numerical imprecision for larger values of $|\xi|$, which is visible in Figure \ref{Fig:num_comp} for the case $\mu = 0.5$.

\section{Results and Discussion}\label{sec.results}

\begin{figure}
\centering
\subfloat[$\mu = 0.5$]{
\includegraphics[width=0.499\textwidth]{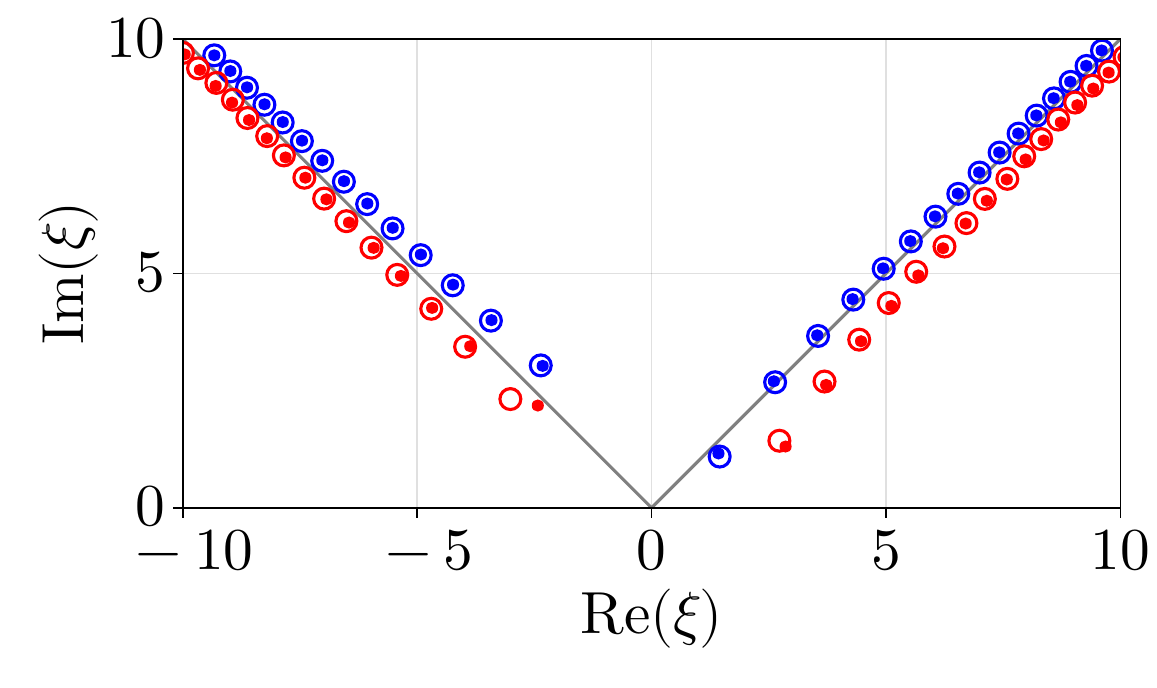}
}
\subfloat[$\mu = 1$]{
\includegraphics[width=0.499\textwidth]{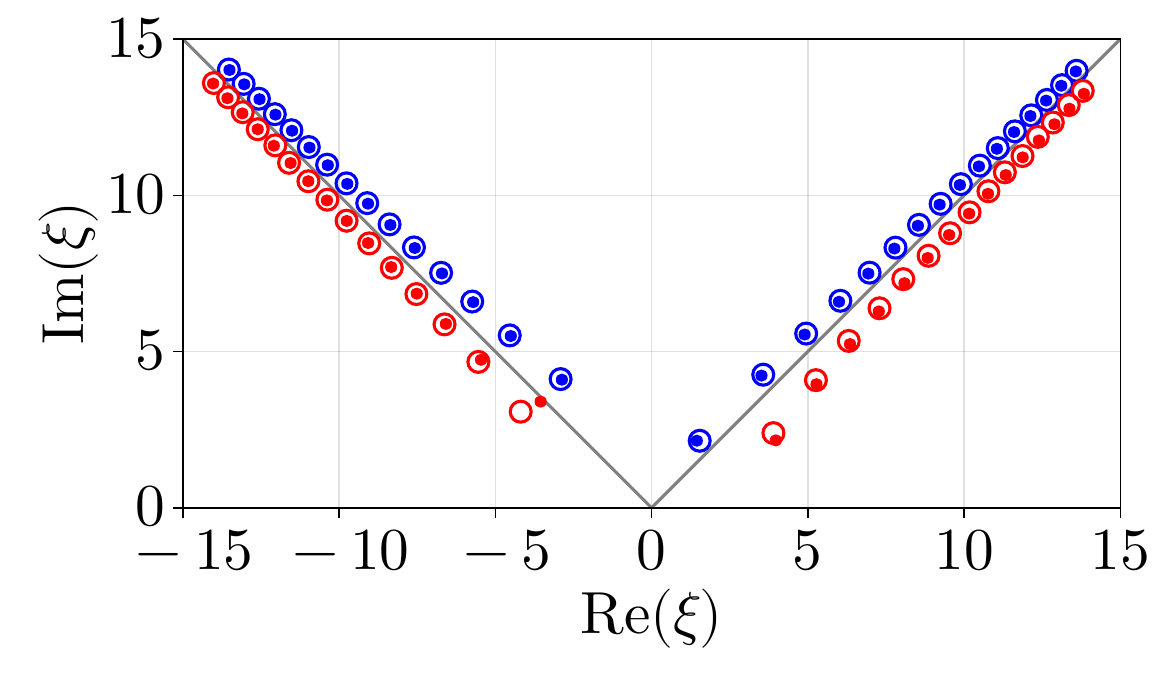}
}

\subfloat[$\mu = 2$]{
\includegraphics[width=0.499\textwidth]{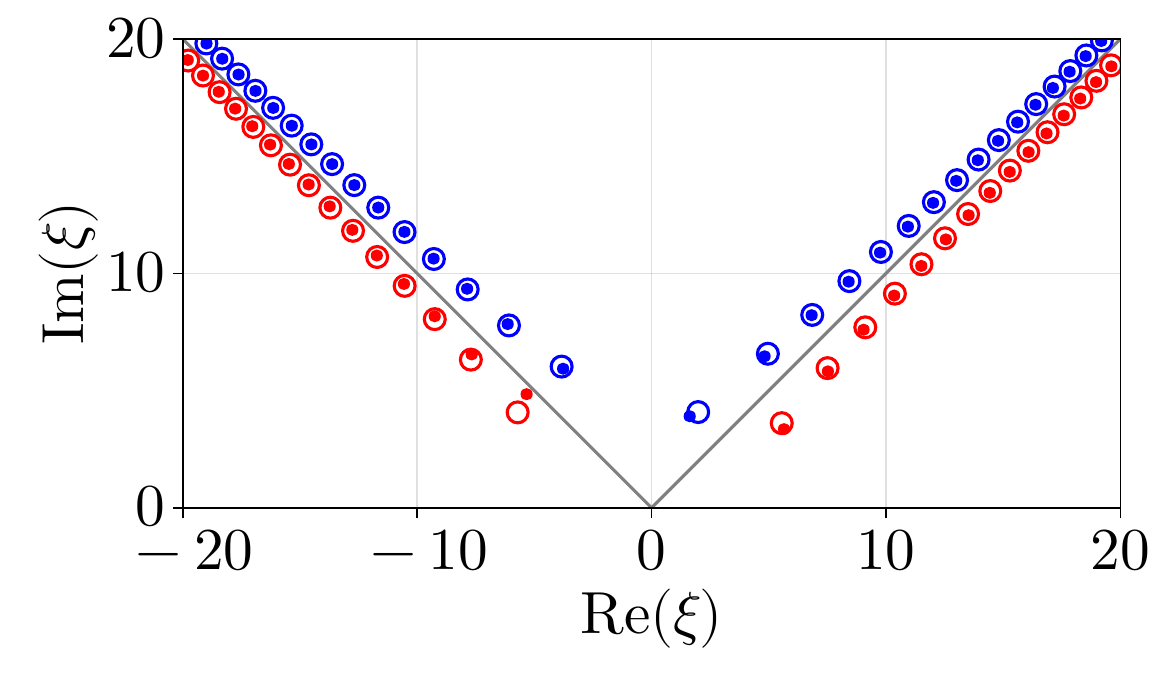}
}
\caption{Comparison of numerical poles and zeroes with asymptotically-calculated pole positions for (a) $\mu = 0.5$, (b) $\mu = 1$, and (c) $\mu = 2$. The numerical poles and zeroes are represented by empty blue and red circles respectively. The asymptotic pole and zero locations from \eqref{e:poleloc} and \eqref{e:zeroloc} are represented by filled blue and red circles respectively. The anti-Stokes curves are shown as gray lines. The asymptotic results provide accurate predictions of the pole and zero locations, with increasing accuracy for larger $|\xi|$.}\label{Fig:num_comp}
\end{figure}

We now compare the transasymptotic predictions \eqref{e:poleloc} and \eqref{e:zeroloc} with the numerically calculated pole and zero locations, obtained numerically from the exact solution (\ref{eq:exactsoln}) given in Ref.~\cite{vandenheuvel2022burgers} (see Section~\ref{sec:exact}).  Examples are provided in Figure \ref{Fig:num_comp} for (a) $\mu = 0.5$, (b) $\mu = 1$, and (c) $\mu = 2$. Note that there is some loss of numerical precision in the calculations for locating zeros via the exact solution for larger values of $|\xi|$, particularly in the $\mu = 0.5$ case. As this is the region in which the asymptotic predictions are most accurate, this provides further motivation for studying poles using both asymptotic and numerical approaches.

For both poles and zeroes, the asymptotic approximation becomes more accurate as $|M|$ increases, corresponding to poles and zeroes that are further from the origin. For smaller values of $|M|$, the pole locations are significantly more accurate than the zero locations; this is expected, as we obtained more terms of accuracy in the approximation for the pole locations \eqref{e:poleloc} than in the approximation for the zero locations \eqref{e:zeroloc}.  Further, for the values of $\mu$ used in Figure~\ref{Fig:num_comp}, it is remarkable that the the approximations \eqref{e:poleloc} and \eqref{e:zeroloc} are particularly effective for identifying the location of the nearest pole and zero to the real axis, which, as we have just noted, are the least accurate.

Returning to Figure \ref{Fig:complex_comp}, we overlay the asymptotic predictions of poles \eqref{e:poleloc} and zeroes \eqref{e:zeroloc} on phase portraits of the exact solution (\ref{eq:exactsoln}) of (\ref{eq.z0ODE})--(\ref{eq.farfield}). We can see the phase changes rapidly in the vicinity of each pole and zero.  Indeed, these plots support the claim made in VandenHeuvel et al.~\cite{vandenheuvel2022burgers} that the simple poles and simple zeroes occur in interlaced pairs. This pattern is explained by the form of \eqref{eq:r-poles} and \eqref{eq:r-zeros}. The poles are offset by $(2M+1)\pi \i$ and the zeroes are offset by $(2M-1/2)\pi \i$, where $M \in \mathbb{Z}$. This will result in each pole being nearest to a particular zero with an offset that differs by $\pi\i/2$.

\begin{figure}
\centering
\subfloat[$\mu = 0.5$]{
\includegraphics[width=0.33\textwidth]{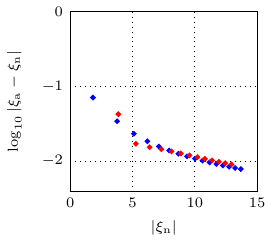}
}
\subfloat[$\mu = 1$]{
\includegraphics[width=0.33\textwidth]{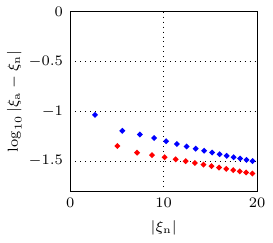}
}
\subfloat[$\mu = 2$]{
\includegraphics[width=0.33\textwidth]{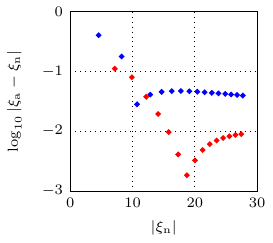}
}
\caption{Comparison of asymptotic and numerical predictions of pole locations (denoted $\xi_{\mathrm{a}}$ and $\xi_{\mathrm{n}}$ respectively), for (a) $\mu = 0.5$, (b) $\mu = 1$, and (c) $\mu = 2$. Poles in the first quadrant are represented by blue circles, while poles in the second quadrant are represented by red circles. In each case, the absolute error decreases as $|\xi_\mathrm{n}|$ increases, due to the increase in $|M|$.}\label{Fig:log_plot}
\end{figure}

In Figure \ref{Fig:log_plot}, we show the error in the pole positions for each value of $\mu$. The location of the poles obtained from the numerical solution is denoted as $\xi_\mathrm{n}$, and the location of the poles obtained from the asymptotic solution is denoted as $\xi_{\mathrm{a}}$. By plotting the difference between the numerical and asymptotic pole locations as a function of $|\xi_\mathrm{n}|$, we see that the accuracy of the approximation increases as $|M|$, and hence $|\xi_{\mathrm{n}}|$ grows. 

Finally, we emphasise that poles and zeros of $U(\xi)$ relate to Burgers' equation (\ref{e:burger}) via $x=\mathrm{i}+t^{1/2}\xi$.  Thus, in the $x$-plane each pole and zero moves with speed $\mathcal{O}(t^{-1/2})$ as $t\rightarrow 0^+$.  For a discussion about the ultimate behaviour of these poles for the initial condition (\ref{e:IC}), see \cite{vandenheuvel2022burgers}. Included in \cite{vandenheuvel2022burgers} is an application of matched asymptotic expansions to determine the leading-order location of poles for $t=\mathcal{O}(1)$ and $t\gg 1$.  This brief analysis exploits the the fact that a series expansion for large $x$ has a leading-order term that is the same size of the exponential correction along anti-Stokes curves, which suggests an appropriate rescaling.  In this way, the matched asymptotic expansions approach in \cite{vandenheuvel2022burgers} for $t=\mathcal{O}(1)$ and $t\gg 1$ is similar in spirit to our own analysis for $t\ll 1$ (although the former only derives a first-order approximation for pole locations).

\section{Conclusions}\label{sec:conclude}

We have studied solutions to (\ref{eq.z0ODE})--(\ref{eq.farfield}), which represents the small-time solution to Burgers' equation \eqref{e:burger} near singular points of the initial condition \eqref{e:IC}. This work was motivated by VandenHeuvel et al.~\cite{vandenheuvel2022burgers}, who used special function theory to show that the analytic continued solution (\ref{eq:exactsoln}) contains two rays of simple poles and simple zeros which rapidly stream out of the singularities of the initial condition \eqref{e:IC} for $t\ll 1$.

We have presented a methodology (based on \cite{costin1999correlation,costin2002formation,costin2015tronquee}) to approximate locations of poles and zeroes in the solution of (\ref{eq.z0ODE})--(\ref{eq.farfield}) for $\mu=\mathcal{O}(1)$, which, in principle, allow us to derive asymptotic expansions to arbitrary asymptotic order. This approach does not rely on the exact solution or special function theory in any way, and can be applied to differential equations which do not have such a convenient  formulation. 

The approach used is as follows:
\begin{itemize}
    \item We applied exponential asymptotics to determine the Stokes structure of the solution, and to calculate the form of exponentially small terms \eqref{eq.Uexp} that appear across Stokes curves.
    \item We expressed the solution as a transseries \eqref{eq.transseries5}, which contains all of the subdominant exponential contributions to the solution. This transseries is asymptotically valid in regions where the exponential terms are small.
    \item We switched the order of summation in the transseries, and summed the expression over all exponential terms at each algebraic order. This allowed us to extend the range of asymptotic validity so that the new series \eqref{eq.transseries4} describes the solution in regions where the exponentials are large. 
\end{itemize}
Using this method, we calculated the location of poles and zeroes in the transseries expression, which we then inverted to obtain asymptotic expressions for the pole locations \eqref{e:poleloc} and zero locations \eqref{e:zeroloc} in terms of the original coordinate $\xi$. These could be extended to include further corrections using additional terms in the transseries (eg. those in \eqref{eq.a13}--\eqref{eq.a02}, or \eqref{e.A2}). 

Another purpose of this study was to illustrate the relationship between exponential asymptotics and transasymptotic analysis, and how they combine to describe the solution to differential equations. This connection is explicit in that the exponentially small contribution described in \eqref{eq.Uexp} was used to specify the transseries parameter $\sigma$. This parameter fixed the transseries coefficients $a_m^{(n)}$, and hence the boundary conditions used to solve the differential equations for $A_0(\tau)$ and $A_1(\tau)$. Exponential asymptotics describes the small exponential contributions that appear across Stokes curves, and transasymptotic analysis allows us to find an asymptotic expression for the solution when these contributions become large across anti-Stokes cruves. Together, the methods provide a rather complete picture of the solution behaviour.

There are extensions to this analysis which were not considered here. As discussed in \cite{vandenheuvel2022burgers}, if the complex singularities in the initial condition for Burgers' equation \eqref{e:burger} are not simple poles, we must introduce an extra boundary layer into the small-time analysis. In general, we expect additional complications of this nature to arise when the strength of the singularity in the initial condition does not match the strength of singularities in the solution of the PDE for $t > 0$ predicted by dominant balance arguments.  For Burgers' equation, this further work will involve the prediction of branch point locations in the solution (see \cite{Aniceto:2022dnm}).

Finally, we note that singularities and zeroes of solutions to ODEs are often arranged in more complicated lattices, rather than rays (as we have been dealing with here for (\ref{eq.z0ODE})). This famously occurs in solutions to the Painlev\'{e} equations (eg. \cite{costin2015tronquee,fornberg_2014}), which frequently describe the solutions to partial differential equations near catastrophes, as in \cite{dubrovin2009universality}.  For solutions to ODEs with a lattice of singularities, the transasymptotic analysis for locating singularities and zeros becomes more complicated; we either have to resolve each ray of the lattice individually \cite{costin2015tronquee}, or use more sophisticated methods to determine the full lattice behaviour (see \cite{Aniceto-upcoming_painleve}).

\section{Acknowledgements}
The authors are grateful for valuable discussions with John R. King. CJL acknowledges the support of Australian Research Council Discovery Project DP190101190. IA acknowledges the support of the EPSRC Early Career Fellowship EP/S004076/1. CJL, IA and SWM would like to thank the Isaac Newton Institute for Mathematical Sciences, Cambridge, for support and hospitality during the programme ``Applicable resurgent asymptotics: towards a universal theory'', where much of the work on this paper was undertaken. The programme was supported by the EPSRC grant no EP/R014604/1.

\appendix

\section{Background on exponential asymptotics and transseries}\label{sec:background}

\subsection{Exponential asymptotics}\label{sec:stokes_background}

Exponential asymptotic methods were developed by Berry~\cite{Berry1988,Berry} and others to study exponentially small contributions to asymptotic expansions for singularly-perturbed problems which Stokes' phenomenon~\cite{Stokes}. This theory describes the rapid appearance of exponentially small contributions as certain curves in the complex plane, known as Stokes curves, are crossed. Exponential asymptotics provides a means to study this switching behaviour in the exponentially small terms.

The version of exponential asymptotics that we use here in our study was developed by Chapman, Olde Daalhuis and co-workers~\cite{Chapman,Daalhuis}. Typically, this method involves writing the solution $V(z; \epsilon)$ of some singularly-perturbed equation as a power series such as (\ref{eq.Vser})
We obtain an expression for $V_n$ by substituting the series~\eqref{eq.Vser} into the governing equation and matching terms at powers of $\epsilon$. To calculate each term $V_n$, we typically differentiate earlier terms in the power series. If the leading-order solution $V_0$ contains singularities in the $z$-plane, the repeated differentiation guarantees that later terms in the series must contain stronger singularities in the same location, and therefore exhibit asymptotic behaviour known as ``factorial-over-power divergence"~\cite{Dingle}. 

To capture the factorial-over-power divergence, Chapman et al.~\cite{Chapman} proposed an asymptotic ansatz for the terms $V_n$ in the limit $n\rightarrow\infty$, known as ``late-order terms", of the form (\ref{e:lateorder_intro}).  
We set $\gamma$ to be constant, and $G$ and $\chi$ are functions of the independent variable but are independent of $n$. The function $\chi$ is known as the ``singulant".  It is equal to $0$ at each singularity of the leading-order solution, meaning that each $V_n$ is also singular at the these points. The function $G$ is known as the ``prefactor''.

By substituting~\eqref{eq.Vser} and~\eqref{e:lateorder_intro} into the governing equation and matching orders of $\epsilon$, we can find $\chi$ and $G$. To determine $\gamma$, we apply the condition that the late-order terms \eqref{e:lateorder_intro} must have the same singularity strength as $V_0$ if we set $n=0$.

The series~\eqref{eq.Vser} is then truncated. A key result from Refs~\cite{Berry1988,Berry} is that if the series is truncated optimally to minimise the truncation error, the remainder will be exponentially small, and can be studied to determine the behaviour of the exponentially small contributions.  We define this truncation point to be $n=N-1$ and obtain 
\begin{equation}
	V = \sum_{n=0}^{N-1} \eps^{2n}V_n+V_{\exp}\,,
\label{e:series_intro_1}
\end{equation}
where $R_N$ is the exponentially small truncation remainder. 

By substituting~(\ref{e:series_intro_1}) into the governing equation, we can obtain a new equation for the exponentially small remainder term. Following~\cite{Daalhuis}, we can apply matched asymptotic expansions in the neighbourhood of the Stokes curve to find that the remainder has the form 
\begin{equation}
	V_{\exp}\sim\mathcal{S}G\e^{-\chi/\epsilon^2} \quad \mathrm{as} \quad \epsilon \rightarrow 0\,,
\label{e:stokes_intro}
\end{equation}
where $\mathcal{S}$ is known as the Stokes multiplier; it is a function of the independent variables in the problem. The Stokes multiplier is essentially constant except in the neighbourhood of Stokes curves, where it undergoes a rapid jump, often switching the contribution on or off entirely. 

From \eqref{e:stokes_intro}, we see that the singulant controls the asymptotic behaviour of the exponential term in the expression. Finding the singulant allows us to determine the location of the Stokes curves, as Stokes switching occurs along the curves where the exponential is most asymptotically subdominant compared to the leading-order algebraic series. This corresponds to the condition
(\ref{e:Stokes_def}).
Stokes curves originate at singular points of $V_0$, where $\chi = 0$, and follow curves on which~\eqref{e:Stokes_def} is satisfied. This method therefore determines both the form of the exponentially small terms \eqref{e:stokes_intro} and the location of the curves across which this behaviour appears \eqref{e:Stokes_def}. 

The other important class of curves is anti-Stokes curves, which satisfy
(\ref{eq.antistokes}).
These are curves where the exponential terms are no longer exponentially small, but are in fact algebraic in the asymptotic limit. The series \eqref{eq.Vser} is not able to accurately describe the solution behaviour beyond these curves, as the series is no longer asymptotically valid. The competition between these growing exponential terms is what will give rise to poles and zeros in the solution of the original problem, so locating these poles and zeros requires continuing the solution past the anti-Stokes curves.  For this purpose, we apply transseries methods.

\subsection{Transseries}

Transasymptotics was introduced in the context of Painlev\'{e} equations to predict the position of moveable poles given boundary data~\cite{costin1999correlation,costin2002formation,costin2015tronquee}. The theory has been effectively used to extend models of relativistic hydrodynamics \cite{Basar:2015ava,Behtash:2020vqk,Aniceto:2022dnm} and field theories \cite{Borinsky:2020vae} beyond their initial regime. Further, it has been used in previous analysis of the Burgers equation for small viscosity \cite{Chapman_shock_caustic}, and to predict the appearance of different scales in discrete bifurcation phenomena \cite{Aniceto_2021}. 

The starting point for our transseries method is constructing the asymptotic series of some solution $U(\xi)$ in powers of a large parameter $\xi$. We denote the algebraic power series for $U(\xi)$ in the limit that $|\xi| \to \infty$ as $U_0$, where
\begin{equation}
	U_0(\xi)\sim \sum_{m=0}^{\infty} \frac{a_m^{(0)}}{\xi^m}.
\end{equation}
The asymptotic expansion of $U(\xi)$ can also contain exponentially small terms, which have asymptotic series of their own that are not contained within $U_0(\xi)$. We represent the exponent of these terms, which will be the singulant from Section \ref{sec:stokes_background}, by $\chi(\xi)$; in some references, this exponent is known as the action. In nonlinear problems, integer multiples of each exponential contribution also appear, with higher integers creating successively smaller exponential terms. We denote the series for each exponential term as $U_{ n}$, where
\begin{equation}
	U_n(\xi)\sim \xi^{ n \beta}\mathrm{e}^{- n \chi(\xi)}\sum_{m=0}^{\infty} \frac{a_m^{( n)}}{\xi^m}.
\end{equation}
If $\mathrm{Re}(\chi) > 0$, the exponential term is small in the asymptotic limit. Typically, exponentially small terms appear as Stokes curves are crossed. If we continue the solution across anti-Stokes curves into regions of the plane where $\mathrm{Re}(\chi) < 0$, the exponential terms become large, and the series terms reorder.

To study the behaviour of the solution in regions beyond anti-Stokes curves, we require an asymptotic description that includes not only its algebraic (or perturbative) power series, but also the exponentially small (or non-perturbative) contributions. We track these contributions by writing an expression for the solution known as a ``transseries'', which is a multiple series expansion in both the large parameter $\xi$ as well as the exponential terms $\mathrm{e}^{-\chi(\xi)}$, or
\begin{equation}
	U(\xi) \sim \sum_{ n=0}^{\infty} \sigma^{ n} U_{n}(\xi) = \sum_{ n=0}^{\infty} \sigma^{ n}\xi^{ n \beta} \mathrm{e}^{- n \chi(\xi)} \sum_{m = 0}^{\infty} \frac{a_m^{( n)}}{\xi^m}.\label{eq:trans-intro}
\end{equation}
We have introduced a transseries parameter $\sigma$ which encodes information about initial or boundary conditions. In regions where there are no exponentially small terms we set $\sigma = 0$. If we continue past a Stokes curve that produces exponentially small contributions, we select a nonzero value of $\sigma$ using the known leading behaviour of the exponential contributions obtained using exponential asymptotics. 

The expression in \eqref{eq:trans-intro} is the simplest transseries that can be written, and is suitable for solutions to nonlinear ordinary differential equations containing multiples of a single exponential. Problems containing multiple exponential contributions, such as \cite{chapmanmortimer_2005}, require transseries that contain integer multiples of each distinct exponential, each of which has a corresponding transseries parameter.

The series associated with each exponential in the transseries \eqref{eq:trans-intro} are all asymptotic, and accurate numerical values can be found using summation techniques  \cite{BerryHowls1990,Berry1991,olde1995hyperasymptotic,olde1995hyperasymptotic2,olde2005hyperasymptotics,Sauzin-summability,ANICETO20191}. The transseries in \eqref{eq:trans-intro} also provides a starting point for analytically continuing past the anti-Stokes curve.

We now use a method known as ``transasymptotic summation'', which corresponds to changing the order of summations in the transseries \eqref{eq:trans-intro} and taking the sum over all exponentials at each power of $\xi^{-1}$. We write
\begin{equation}\label{eq:transasymp-intro}
	U(\xi) \sim \sum_{m = 0}^{\infty} \frac{1}{\xi^m} \sum_{ n = 0}^{\infty} 
	\left(\sigma \xi^{\beta} \mathrm{e}^{-\chi(\xi)}\right)^{ n} a_m^{( n)} = \sum_{m = 0}^{\infty} \frac{A_m(\tau)}{\xi^m},
\end{equation}
where we introduced the new variable $\tau$ and functions $A_m$ such that
 \begin{equation}
 	\tau=\sigma\xi^{\beta}\mathrm{e}^{-\chi(\xi)}, \qquad A_m(\tau)=\sum_{ n=0}^{\infty} \tau^ n a_m^{( n)}.\label{eq:tau-def-intro}
\end{equation}
Note that this equality for $A_m$ is only valid in sectors where $\tau \to 0$ as $|\xi| \to \infty$; outside of these sectors, $A_m$ is equal to the analytic continuation of the summed expression. This definition allows us to continue the summed expression \eqref{eq:transasymp-intro} beyond anti-Stokes curves and into regions of the plane where $\tau$ is not necessarily small in the asymptotic limit. Hence, while the original transseries was obtained with the assumptions $\mathrm{e}^{-\chi(\xi)}\ll \,|\xi|^{-1}\ll 1$, the transasymptotic summation allows us to study regions where the parameter $\xi$ is large\footnote{Although these approximations are typically accurate even for values of $|\xi|$ that are not particularly large.}, but the exponentials are not small. 
 
 The functions $A_m(\tau)$ satisfy ordinary differential equations with boundary conditions that depend on $a_m^{( n)}$. In general, poles and zeros in $\tau$ will be introduced by the nonlinear differential equation for $A_0$. Subsequent terms $A_m$ for $m \geq 1$ are generated by linear equations and do not introduce new poles; instead, these terms provide asymptotic corrections to the locations of the poles that appear in the leading-order $A_0$. 
 
 Finally, we will express the pole and zero locations in $\tau$ as a series in $\xi^{-1}$, and invert the expression to obtain an asymptotic description of the pole locations in terms of $\xi$ that is valid as $|\xi| \to \infty$. We will apply the same general procedure to also determine the location of zeroes in the solution.

\section{Calculating \texorpdfstring{$\Lambda$}{lambda}}\label{s.lambda}

\begin{figure}
\centering
\subfloat[Estimate of $\mathrm{Re}(\Lambda)$ after $n$ iterations.]{
\includegraphics[]{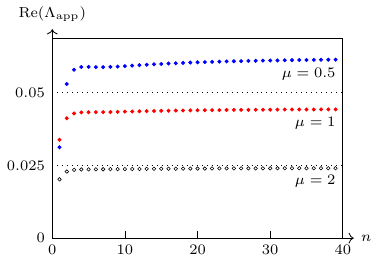}
}\quad
\subfloat[Estimate of $\mathrm{Im}(\Lambda)$ after $n$ iterations.]{
\includegraphics[]{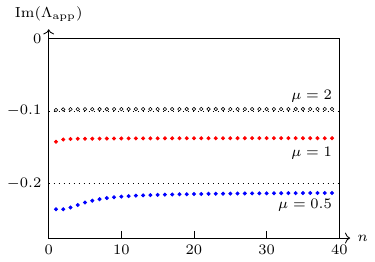}
}
\caption{The approximate value of $\Lambda$, denoted as $\Lambda_{\mathrm{app}}$, obtained using the value of $\hat{V}_n$ over a range of $n$ in the expression \eqref{eq.Vcomp}. In each case, $\Lambda_{\mathrm{app}}$ converges to a constant value as $n$ grows. }\label{Fig.Lambda1}
\end{figure}

In order to determine the late-order terms $V_n$ in (\ref{e:lateorder_intro}), we require the value of $\Lambda$, or the constant associated with the prefactor. To accomplish this, we note that the series expression \eqref{eq.Vser} fails to be asymptotic when consecutive terms are of equal size in the limit $\eps \to 0$. From the late-order terms  (\ref{e:lateorder_intro}), we see that this occurs if $z^2/4\mu = \mathcal{O}(\epsilon^2)$ in this limit. This suggests that a sensible re-scaling is given by $z/2\sqrt{\mu} = \epsilon \hat{z}$, which we apply to \eqref{eq.z0ODE}. We define a new variable governing the solution behaviour in this region, $V(z) = \hat{V}(\hat{z})/\epsilon$, giving the re-scaled equation 
\begin{equation}\label{eq.Vinn}
\frac{1}{2}\diff{^2 \hat{V}}{\hat{z}^2} - \frac{1}{\sqrt{\mu}} \hat{V} \diff{\hat{V}}{\hat{z}}  + \hat{z}\diff{\hat{V}}{\hat{z}}+ \hat{V}= 0.
\end{equation}
To match with the late-order terms, we must determine the behaviour of the solution in this inner region in the limit $|\hat{z}|\to\infty$. We therefore pose the series
\begin{equation}\label{eq.Vinsum}
\hat{V}(\hat{z}) \sim \sum_{n=0}^{\infty}\frac{\hat{V}_n}{\hat{z}^{2n+1}},\quad \mathrm{as} \quad |z| \to \infty.
\end{equation}
Applying the series \eqref{eq.Vinsum} to the differential equation \eqref{eq.Vinn} gives $\hat{V}_0$ being arbitrary, and 
\begin{equation}
\hat{V}_n = \left(n-\frac{1}{2}\right)\hat{V}_{n-1} + \frac{1}{\sqrt{\mu}} \sum_{j=0}^{n-1}\left(\frac{j + 1/2}{n}\right) \hat{V}_j \hat{V}_{n-j-1}, \qquad n \geq 1.
\end{equation}
By rewriting the inner expansion \eqref{eq.Vinsum} in terms of the outer variable $z$ and comparing it with the late-order term ansatz  (\ref{e:lateorder_intro}), we see that the matching condition is given by 
\begin{equation}\label{eq.Vcomp}
\Lambda = \lim_{n \to \infty}\frac{\hat{V}_n (4\mu)^{\i/4\mu}}{\Gamma\left(n + \tfrac{1}{2}-\tfrac{\i}{4\mu}\right)}.
\end{equation}
To determine $\hat{V}_0$, we note that $\hat{V}_n$ must grow factorially in $n$ as $n \to \infty$ in order for \eqref{eq.Vcomp} to converge. Assuming that $n$ grows in this fashion, we see that for large $n$,
\begin{equation}
\hat{V}_n \sim \left(n - \frac{1}{2} + \frac{\hat{V}_0}{\sqrt{\mu}}\right) \hat{V}_{n-1} \quad \mathrm{as} \quad n \to \infty.
\end{equation}
Comparing this to \eqref{eq.Vcomp} gives the requirement that $\hat{V}_0 = -\i/4\sqrt{\mu}$, in order for the factorial growth to be consistent between the inner limit of the outer expansion and the outer limit of the inner expansions. We can now use \eqref{eq.Vcomp} to determine numerical values of $\Lambda$ by calculating $\hat{V}_n$ for large values of $n$, and testing to see the value to which the ratio converges. In Figure \ref{Fig.Lambda1}, we illustrate this approximation procedure for $\mu = 0.5$, 1, and 2. In Figure \ref{Fig.Lambda2}, we illustrate values of $\Lambda$ over a range of $\mu$ values, obtained by evaluating the ratio in \eqref{eq.Vcomp} at $n = 1000$ in order to approximate $\Lambda$.

\bibliography{reference2} 
\bibliographystyle{plain}

\end{document}